%%%% prevent double loading:
\expandafter\ifx\csname mthreemacsloaded\endcsname\relax\else \fi

\magnification1100
\input amstex

%%% Hack of Plain TeX correction and style macros 
%%% written by Walter Neumann and Larry Siebenmann:

 \catcode`\@=11
 \let\wlog@ld\wlog
 \def\wlog#1{\relax}

 \newif\ifIN@
 \def\m@rker{\m@@rker}
 \def\IN@{\expandafter\INN@\expandafter}
 \long\def\INN@0#1@#2@{\long\def\NI@##1#1##2##3\ENDNI@
    {\ifx\m@rker##2\IN@false\else\IN@true\fi}%
     \expandafter\NI@#2@@#1\m@rker\ENDNI@}
  \newtoks\Initialtoks@  \newtoks\Terminaltoks@
  \def\SPLIT@{\expandafter\SPLITT@\expandafter}
  \def\SPLITT@0#1@#2@{\def\TTILPS@##1#1##2@{%
     \Initialtoks@{##1}\Terminaltoks@{##2}}\expandafter\TTILPS@#2@}
  \newtoks\Trimtoks@

 \def\ForeTrim@{\expandafter\ForeTrim@@\expandafter}
 \def\ForePrim@0 #1@{\Trimtoks@{#1}}
 \def\ForeTrim@@0#1@{\IN@0\m@rker. @\m@rker.#1@%
     \ifIN@\ForePrim@0#1@%
     \else\Trimtoks@\expandafter{#1}\fi}
 
  \def\Trim@0#1@{%
      \ForeTrim@0#1@%
      \IN@0 @\the\Trimtoks@ @%
        \ifIN@
             \SPLIT@0 @\the\Trimtoks@ @\Trimtoks@\Initialtoks@
             \IN@0\the\Terminaltoks@ @ @%
                 \ifIN@
                 \else \Trimtoks@ {FigNameWithSpace}%
                 \fi
        \fi
      }

  %%% Math Bolds
  \font\titlebold=cmbx12 scaled 1200
  \font\twelvebold=cmbx12
  \font\tenbold=cmbx10
  \font\ninebold=cmbx9
  \font\sevenbold=cmbx7
  \font\fivebold=cmbx5

  \input amssym.def \input amssym
  %%% point sizes not loaded by amssym.def:
     \font\titlemsa=msam10 at 14.4pt
     \font\titlemsb=msbm10 at 14.4pt
     \font\titleeufm=eufm10 at 14.4pt
     \font\twelvemsa=msam10 scaled 1200
     \font\twelvemsb=msbm10 scaled 1200
     \font\twelveeufm=eufm10 scaled 1200
     \font\ninemsa=msam9
     \font\ninemsb=msbm9
     \font\nineeufm=eufm9

   %%% Cyrillic fonts (for accents and input, see ams cyr doc)
   \ifx\cyrfam\undefined
   \else
     \immediate\write16{}%
     \message{ !!! cyr fonts already defined. !!! }
     \message{ --- edit out superfluous font defs? }
   \fi
   \newfam\cyrfam
       \font\titlecyr=wncyr10 scaled 1440 %%% no caps?
       \font\twelvecyr=wncyr10 scaled 1200
       \font\tencyr=wncyr10
       \font\ninecyr=wncyr9
       \font\sevencyr=wncyr7
       \font\sixcyr=wncyr6

   %%% Euler script fonts (replacing caligraphic):
   \newfam\eusmfam
       \font\titleeusm=eusm10 scaled 1440
       \font\twelveeusm=eusm10 scaled 1200
       \font\teneusm=eusm10
       \font\nineeusm=eusm9
       \font\seveneusm=eusm7
       
       \font\fiveeusm=eusm5

\let\Cal\cal

 %%% Some fonts not loaded by plain
    \font\ninemrm=cmr9 %% new name for 9 pt math roman
    \font\ninei=cmmi9
    \font\ninesy=cmsy9 
    \skewchar\ninei='177
    \skewchar\ninesy='60

  \font\twelvemrm=cmr10 at 12pt %% new name
  \font\twelvei=cmmi10 at 12pt
  \font\twelvesy=cmsy10 at 12pt
 % \font\twelveex=cmex10 at 12pt

  \font\titlemrm=cmr10 at 14.4pt %% new name
  \font\titlei=cmmi10 at 14.4pt
  \font\titlesy=cmsy10 at 14.4pt
 % \font\titleex=cmex10 at 14.4pt

 %%%% Miscellanious font definitions

  \def\Smallfonts{\ninepoint}

  \def\Hfont{\titlepoint\bf}
  \def\Authorfont{\twelvepoint\it}
  \def\HHfont{\twelvepoint\bf}
  \def\HHHfont{\bf}
  % automatically smaller in 9 point parts
  \def\Bibfont{\tenbf}
  \def\Coordfont{\nineit }% defined in osuPSfnt.sty

  \def \thfont {\bf }
  \def \pffont {\it\itSpacing }
  \def \rkfont {\bf }
  \def \dffont {\bf }
  \def \egfont {\bf }

 %%%%% NINEPOINT %%%%%
 \def\ninepoint{%
  \def\rm{\fam0\ninerm}%
    \textfont0=\ninemrm  \scriptfont0=\sevenrm  \scriptscriptfont0=\fiverm
    \textfont1=\ninei    \scriptfont1=\seveni   \scriptscriptfont1=\fivei
  \def\mit{\fam1\ninei}%
  \def\oldstyle{\fam1\ninei}%
    \textfont2=\ninesy   \scriptfont2=\sevensy  \scriptscriptfont2=\fivesy
    \textfont3=\tenex    \scriptfont3=\tenex    \scriptscriptfont3=\tenex
  \def\it{\fam\itfam\nineit}%
    \textfont\itfam=\nineit
  \def\bf{\ifmmode\fam\bffam\else\ninebf\fi}%
    \textfont\bffam=\ninebold 
    \scriptfont\bffam=\sevenbold 
    \scriptscriptfont\bffam=\fivebold%
  \def\msa{\fam\msafam\ninemsa}%
    \textfont\msafam=\ninemsa 
    \scriptfont\msafam=\sevenmsa
    \scriptscriptfont\msafam=\fivemsa%
  \def\msb{\fam\msbfam\ninemsb}%
    \textfont\msbfam=\ninemsb%
    \scriptfont\msbfam=\sevenmsb%
    \scriptscriptfont\msbfam=\fivemsb%
  \def\eufm{\fam\eufmfam\nineeufm}%
    \textfont\eufmfam=\nineeufm
    \scriptfont\eufmfam=\seveneufm
    \scriptscriptfont\eufmfam=\fiveeufm
   \def\eusm{\fam\eusmfam\nineeusm}%
     \textfont\eusmfam=\nineeusm
     \scriptfont\eusmfam=\seveneusm
     \scriptscriptfont\eusmfam=\fiveeusm
   \def\cyr{\fam\cyrfam\ninecyr}%
     \textfont\cyrfam=\ninecyr
     \scriptfont\cyrfam=\sevencyr
     \scriptscriptfont\cyrfam=\sixcyr%%
  \setbox\strutbox=\hbox{\vrule
      height7pt depth3pt width0pt}%
   \baselineskip=10.8pt\rm}

 \let\eightpoint\ninepoint % we do not use eightpoint

 %%%%% FONTS AT TENPOINT %%%%%
 \def\tenpoint{%
  \def\rm{\fam0\tenrm}%
    \textfont0=\tenmrm \scriptfont0=\sevenrm \scriptscriptfont0=\fiverm%
  \def\mit{\fam1\teni}%
  \def\oldstyle{\fam1\teni}%
    \textfont1=\teni   \scriptfont1=\seveni  \scriptscriptfont1=\fivei%
    \textfont2=\tensy  \scriptfont2=\sevensy \scriptscriptfont2=\fivesy%
    \textfont3=\tenex  \scriptfont3=\tenex   \scriptscriptfont3=\tenex%
  \def\it{\fam\itfam\tenit}%
    \textfont\itfam=\tenit%
  \def\bf{\ifmmode\fam\bffam\else\tenbf\fi}%
    \textfont\bffam=\tenbold% was tenbold for osu
    \scriptfont\bffam=\sevenbold%
    \scriptscriptfont\bffam=\fivebold%
  \def\msa{\fam\msafam\tenmsa}%
    \textfont\msafam=\tenmsa%
    \scriptfont\msafam=\sevenmsa%
    \scriptscriptfont\msafam=\fivemsa%
  \def\msb{\fam\msbfam\tenmsb}%
    \textfont\msbfam=\tenmsb%
    \scriptfont\msbfam=\sevenmsb%
    \scriptscriptfont\msbfam=\fivemsb%
  \def\eufm{\fam\eufmfam\teneufm}%
   \textfont\eufmfam=\teneufm
   \scriptfont\eufmfam=\seveneufm
   \scriptscriptfont\eufmfam=\fiveeufm
   \def\eusm{\fam\eusmfam\teneusm}%
    \textfont\eusmfam=\teneusm
    \scriptfont\eusmfam=\seveneusm
    \scriptscriptfont\eusmfam=\fiveeusm
   \def\cyr{\fam\cyrfam\tencyr}%
    \textfont\cyrfam=\tencyr
    \scriptfont\cyrfam=\sevencyr
    \scriptscriptfont\cyrfam=\sixcyr%%
  \setbox\strutbox=\hbox{\vrule %
      height8.5pt depth3.5ptwidth0pt}%
  \baselineskip=\StdBaselineskip\rm}

 %%%%% FONTS AT TWELVEPOINT %%%%%
 \def\twelvepoint{%
  \def\rm{\fam0\twelverm}%
    \textfont0=\twelvemrm \scriptfont0=\tenmrm \scriptscriptfont0=\sevenrm
    \textfont1=\twelvei   \scriptfont1=\teni   \scriptscriptfont1=\seveni
  \def\mit{\fam1\twelvei}%
  \def\oldstyle{\fam1\twelvei}%
    \textfont2=\twelvesy  \scriptfont2=\tensy  \scriptscriptfont2=\sevensy
    \textfont3=\tenex  \scriptfont3=\tenex  \scriptscriptfont3=\tenex
  \def\it{\fam\itfam\twelveit}%
    \textfont\itfam=\twelveit
  \def\bf{\ifmmode\fam\bffam\else\twelvebf\fi}%
    \textfont\bffam=\twelvebold
    \scriptfont\bffam=\tenbold%
    \scriptscriptfont\bffam=\sevenbold%
  \def\msa{\fam\msafam\twelvemsa}%
    \textfont\msafam=\twelvemsa%
    \scriptfont\msafam=\tenmsa%
    \scriptscriptfont\msafam=\sevenmsa%
  \def\msb{\fam\msbfam\twelvemsb}%
    \textfont\msbfam=\twelvemsb%
    \scriptfont\msbfam=\tenmsb%
    \scriptscriptfont\msbfam=\sevenmsb%
  \def\eufm{\fam\eufmfam\twelveeufm}%
   \textfont\eufmfam=\twelveeufm
   \scriptfont\eufmfam=\teneufm
   \scriptscriptfont\eufmfam=\seveneufm
   \def\eusm{\fam\eusmfam\twelveeusm}%
    \textfont\eusmfam=\twelveeusm
    \scriptfont\eusmfam=\teneusm
    \scriptscriptfont\eusmfam=\seveneusm
   \def\cyr{\fam\cyrfam\tencyr}%
    \textfont\cyrfam=\twelvecyr
    \scriptfont\cyrfam=\tencyr
    \scriptscriptfont\cyrfam=\sevencyr%%
  \setbox\strutbox=\hbox{\vrule
      height10.2pt depth4.55pt width0pt}%
  \baselineskip=14pt\rm}

 %%%%% FONTS AT TITLEPOINT %%%%%
 \def\titlepoint{%
    \textfont0=\titlemrm \scriptfont0=\twelvemrm \scriptscriptfont0=\tenmrm
    \textfont1=\titlei   \scriptfont1=\twelvei   \scriptscriptfont1=\teni
  \def\mit{\fam1\titlei}%
  \def\oldstyle{\fam1\titlei}%
    \textfont2=\titlesy  \scriptfont2=\twelvesy  \scriptscriptfont2=\tensy
    \textfont3=\tenex% math ext not avail in varying sizes??
    \scriptfont3=\tenex
    \scriptscriptfont3=\tenex
  \def\it{\fam\itfam\titleit}%
    \textfont\itfam=\titleit
  \def\bf{\ifmmode\fam\bffam\else\titlebf\fi}%
    \textfont\bffam=\titlebold
    \scriptfont\bffam=\twelvebold%
    \scriptscriptfont\bffam=\tenbold%
  \def\msa{\fam\msafam\titlemsa}%
    \textfont\msafam=\titlemsa%
    \scriptfont\msafam=\twelvemsa%
    \scriptscriptfont\msafam=\tenmsa%
  \def\msb{\fam\msbfam\titlemsb}%
    \textfont\msbfam=\titlemsb%
    \scriptfont\msbfam=\twelvemsb%
    \scriptscriptfont\msbfam=\tenmsb%
  \def\eufm{\fam\eufmfam\titleeufm}%
    \textfont\eufmfam=\titleeufm
    \scriptfont\eufmfam=\twelveeufm
    \scriptscriptfont\eufmfam=\teneufm
   \def\eusm{\fam\eusmfam\titleeusm}%
     \textfont\eusmfam=\titleeusm
     \scriptfont\eusmfam=\twelveeusm
     \scriptscriptfont\eusmfam=\teneusm
   \def\cyr{\fam\cyrfam\tencyr}%
    \textfont\cyrfam=\titlecyr
    \scriptfont\cyrfam=\twelvecyr
    \scriptscriptfont\cyrfam=\tencyr%%
  \setbox\strutbox=\hbox{\vrule
      height12.3pt depth5.54pt width0pt}%
  \baselineskip=16pt\rm}

 %%%% RUNNING HEADINGS
\newbox\AuthorBox\newbox\TitleBox
\newbox\TFLinebox
\newbox\FLinebox
\newbox\HLinebox
\def\SetTFLinebox#1{\setbox\TFLinebox=\hbox{#1}}
\def\SetFLinebox#1{\setbox\FLinebox=\hbox{#1}}
\def\SetHLinebox#1{\setbox\HLinebox=\hbox{#1}}

 \def\SetAuthorHead#1{%
     \setbox\AuthorBox=\hbox{\ninepoint \it 
           \ignorespaces\frenchspacing#1\unskip}}
 \def\SetTitleHead#1{%
     \setbox\TitleBox=\hbox{\ninepoint \it
           \ignorespaces\frenchspacing#1\unskip}}

 %% Italic Spacing Correction
  \def\itSpacing{\relax}
  \def\itSpacingOff{\relax}

  %% Main section headings

 \def\Hrule{\hrule width0pt height0pt}

 %% skip used around proclamations, after section headings,
  % and before subsection-headings:
  \newskip\ProcSkip \ProcSkip 8pt plus2pt minus2pt

 \newskip\LastSkip
 \def\SaveLastSkip{\LastSkip\lastskip}
 \def\RestoreLastSkip{\vskip-\LastSkip\vskip\LastSkip}

 %% Do not indent next paragraph after a header:
 \def\NoindentAfter{\everypar={\setbox0=\lastbox\everypar={}}}

 \long\def\H#1\par#2\par{\notenumber=0 \titlepagetrue%
    {
    \baselineskip=20pt
    \parindent=0pt\parskip=0pt\frenchspacing
    \leftskip=0pt plus .2\hsize minus .3\hsize
    \rightskip=0pt plus .2\hsize minus .3\hsize
 \def\\{\unskip\break}%
    \pretolerance=10000 \Hfont #1\unskip\break
     \vskip7pt\Hrule
\hfill \Authorfont #2\hfill\hfill\unskip}
    \vskip48pt plus 4pt minus 4pt% 60pt=48+12pt
    \par\NoindentAfter\rm}

 \long\def\Hi#1\par#2\par{\notenumber=0 \titlepagetrue%
    {  \baselineskip=0pt  \parindent=0pt\parskip=0pt\frenchspacing
    \leftskip=0pt plus .2\hsize minus .3\hsize
    \rightskip=0pt plus .2\hsize minus .3\hsize
}
    \rm}

 %%% Minor section headings

 \newdimen\PageRemainder
  \def\SetPageRemainder{%\maxdimen case at page tops 12-91 LS
     \PageRemainder=\pagegoal
     \ifdim\PageRemainder=\maxdimen\PageRemainder=\vsize
     \else\advance\PageRemainder by -1\pagetotal\fi}

  \def\Rpt@{}\def\Rpt@@{}

  \long\def\HH#1\par{\par%A
  \SaveLastSkip\removelastskip\goodbreak
  \ifdim\LastSkip<30pt %24pt
     \LastSkip 30pt%24pt 
plus 3pt minus 2pt\fi
  \SetPageRemainder\advance\PageRemainder-\LastSkip
  \ifdim\PageRemainder<150pt
       \edef\Rpt@{remain = \the\PageRemainder\noexpand\\
                pagetotal=\the\pagetotal\noexpand\\
                           pagegoal=\the\pagegoal}%
          \fi
   \ifdim\PageRemainder<65pt %%Head plus 4 lines (approx)
       \ifdim\PageRemainder > 0pt
          \edef\Rpt@@{\noexpand\\
                      Had HH PageRemainder$<$\relax 65pt\noexpand\\
                      Hence forced break!}%
     \vskip 0pt plus .2\PageRemainder\eject %% Pull it out a bit
    \fi\fi
    \vskip\LastSkip\Hrule %%%%%%%%\Hrule added
    \pretolerance=10000\rightskip=0pt plus 3em%B
    \hangafter1 \hangindent=2.2em%
    \noindent
    \HHfont \unskip \Ednote{\Rpt@\Rpt@@}%
            \def\Rpt@{}\def\Rpt@@{}%
            \ignorespaces
            #1\par\rightskip=0pt\pretolerance=\StdPretolerance%
    \NoindentAfter
\tenpoint\rm%
     \medskip \vskip\ProcSkip}%interlineskip adds 2pt to this

  \long\def\HHH#1\par{\par%
  \SaveLastSkip\removelastskip\goodbreak
  \ifdim\LastSkip<\ProcSkip%
     \LastSkip\ProcSkip\fi
  \SetPageRemainder\advance\PageRemainder-\LastSkip
  \ifdim\PageRemainder<150pt
       \edef\Rpt@{remain = \the\PageRemainder\noexpand\\
                pagetotal=\the\pagetotal\noexpand\\
                           pagegoal=\the\pagegoal}%
       \fi
   \ifdim\PageRemainder<48pt  %% 4 lines
        \ifdim\PageRemainder > 0pt
             \edef\Rpt@@{\noexpand\\
                      Had HHH PageRemainder$<$\relax48pt\noexpand\\
                      Hence forced break!}%
       \vskip 0pt plus .2\PageRemainder\eject %% Pull it out a bit
      \fi\fi
   \vskip\LastSkip\par\noindent
   \HHHfont \unskip\Ednote{\Rpt@\Rpt@@}%
  \def\Rpt@{}\def\Rpt@@{}%
  \ignorespaces
   #1\unskip.\quad\rm\ignorespaces
   \ignorepars}

  \long\def\ignorepars#1\par{\def\Test{#1}%
     \ifx\Test\Empty\def\This{\ignorepars}%
        \else\def\This{\Test\par}\fi
           \This}
  \def\Empty{}

 \def\Abstract#1\par{\bgroup\Smallfonts\narrower\HHH #1\par}
 \def\endAbstract{\par\egroup}

 %%%%% Proclamations %%%%%

 \def\ProcBreak{\par%
    \ifdim\lastskip<8pt%
    \removelastskip%
    \penalty-200\vskip\ProcSkip\fi}

 \def\th#1\par{\ProcBreak \noindent
   {\thfont\ignorespaces
    #1\unskip.}\it\itSpacing\kern.4em\ignorepars}%\everymath{\/}

 \def\endth{\ProcBreak\rm\itSpacingOff }%\everymath{}

  %% the theorem statement will be in italic by default

 \def\pf#1\par{\ProcBreak %
    \noindent\pffont#1\unskip.\rm\itSpacingOff{\kern .7em}\ignorepars}

 \def\endpf{\medskip \ProcBreak } %% \qed is alternative

  %% A Box for the QED
  \def\qedbox{\hbox{\vbox{
    \hrule width0.2cm height0.2pt
    \hbox to 0.2cm{\vrule height 0.2cm width 0.2pt
             \hfil\vrule height0.2cm width 0.2pt}
    \hrule width0.2cm height 0.2pt}\kern1pt}}

  %% Typing in \qed makes the qedbox right justified:
  \def\qed{\ifmmode\qedbox
    \else\unskip\ \hglue0mm\hfill\qedbox\ProcBreak\fi}

  \def \rk #1\par{\ProcBreak
     \noindent{\rkfont\ignorespaces #1\unskip.}%
     \rm\kern.6em\ignorepars}

  \def \df #1\par{\ProcBreak
     \noindent{\dffont\unskip\ignorespaces #1\unskip.}%
     \rm\kern.6em\ignorepars}

  \def \enddf {\medskip\ProcBreak }

  \def \eg #1\par{\ProcBreak
     \noindent\egfont\unskip\ignorespaces #1\unskip.
     \rm\kern.6em\ignorepars}

  \def \endeg {\medskip\ProcBreak }

  \newdimen\Overhang

   \def\MaxTag@#1#2#3#4#5{\setbox0=\hbox{#4\ignorespaces#2\unskip}%
     \dimen0=\wd0\advance\dimen0 by#3
     \ifdim\dimen0<#5\relax\dimen0=#5\fi
     \expandafter\edef\csname #1Hang\endcsname{\the\dimen0}}

 \def\MaxItemTag#1{\MaxTag@{Item}{#1}{.4em}{\ItemStyle}{\parindent}}%
 \def\MaxItemItemTag#1{%
        \MaxTag@{ItemItem}{#1}{.4em}{\ItemItemStyle}{\parindent}}
 \def\MaxNrTag#1{\MaxTag@{Nr}{#1}{.5em}{\NrStyle}{\parindent}}
 \def\MaxReferenceTag#1{%
        \MaxTag@{Reference}{[#1]}{.6em}{\ninerm}{\parindent}}
 \def\MaxFootTag#1{\MaxTag@{Foot}{#1}{.4em}{\ninerm}{\z@}}

  %% \SetOverhang@ will prevent for tag-text collision
  \def\SetOverhang@{\Overhang=.8\dimen0%
     \advance\Overhang by \wd0\relax%nec!
     \ifdim\Overhang>\hangindent\relax%nec!
       \advance\Overhang by .25\dimen0%
       \Ednote{Tag is pushing text.}\osumess{Tag is pushing text.}%
     \else\Overhang=\hangindent
     \fi}

   %%% \Item
   \def\Item#1{\par\noindent
      \hangafter1\hangindent=\ItemHang
      \setbox0=\hbox{\ItemStyle\ignorespaces#1\unskip}%
      \dimen0=.4em\SetOverhang@% dimen0 is extra space
      \rlap{\box0}\kern\Overhang\ignorespaces}

   %%% \ItemItem
   \def\ItemItem#1{\par\noindent
      \hangafter1\hangindent=\ItemItemHang
      \setbox0=\hbox{\ItemItemStyle\ignorespaces#1\unskip}%
      \dimen0=.4em\SetOverhang@
      \advance\hangindent by \ItemHang
      \kern\ItemHang\rlap{\box0}%
      \kern\Overhang\ignorespaces}

  %%%% \Nr Items without hanging indentation
  \def\Nr#1{\par\noindent\hangindent=\NrHang %not really a hang
    \setbox0=\hbox{\NrStyle\ignorespaces#1\unskip}%
    \dimen0=.5em\SetOverhang@% dimen0 is extra space
    \rlap{\box0}\kern\Overhang
    \hangindent=\z@\ignorespaces}

  %%%% Roster (not compulsory)
  %%  endRoster has to remember \lastskip (e.g. from a \qed) through \egroup.
   \newskip\Rosterskip\Rosterskip 1pt plus1pt %% modifiable
   \def\Roster{\par\ifdim\lastskip<\Rosterskip\removelastskip\vskip\Rosterskip\fi
    \bgroup}
   \def\endRoster{\par\global\edef\LastSkip@{\the\lastskip}\removelastskip
       \egroup\penalty-50\LastSkip\LastSkip@\relax
       \ifdim\LastSkip<\Rosterskip\LastSkip\Rosterskip\fi
       \vskip\LastSkip}%%changed Feb/5/92 WN

 %%%%% Emphasis %%%%%

 %%%%% Vertical spacing %%%%%

 %%%%% References %%%%%

 \def\cite#1{%\relaxnext@
    \def\nextiii@##1,##2\end@{{\frenchspacing\rm 
      \lBr\ignorespaces##1\unskip{\rm,~\ignorespaces##2}\rBr}}%
    \IN@0,@#1@%
    \ifIN@\def\next{\nextiii@#1\end@}\else
    \def\next{{\rm\lBr#1\rBr}}\fi\next}

 %%%%% Bibliography %%%%%

   \def \Bib#1\par{%
       \par\removelastskip\SetPageRemainder
       \ifdim\PageRemainder < 97pt
        \ifdim\PageRemainder > 0pt
        \vfill\eject
       \fi\fi
    \ProcBreak \par\begingroup\parskip=0 pt%
    \goodbreak \vskip 15 pt plus 10 pt
    \noindent\null\hfill\Bibfont% \kern??pt%  (center over what?)
      \ignorespaces #1\unskip\hfill\null\par 
    \frenchspacing \Smallfonts\rm
    \parskip=2.5 pt plus 1 pt minus.5pt%
    \nobreak\vskip 12pt plus 2pt minus2pt\nobreak
    \leftskip=0 pt \baselineskip=10.5pt}

 \def\ReferenceTagSlide{0em}
  \def\ReferenceTagGap{.5em}

  \def \rf#1{\par\noindent
     \hangafter1\hangindent=\ReferenceHang      
     \setbox0=\hbox{\ninerm[\ignorespaces#1\unskip]}%        
     \dimen0=\ReferenceTagGap\SetOverhang@
     \rlap{\kern\ReferenceTagSlide\box0}%       
     \kern\Overhang\ignorespaces}

  \def\ref#1\par#2\par#3\par#4\par{%
     \rf{#1}#2\unskip,\ #3\unskip,\
     #4\unskip.}

  \def\endBib{\par\endgroup\vskip 12pt minus 6pt }

 %%%%% Coordinates %%%%%

  \long\def\Coordinates#1\endCoordinates{%\relax}
 {\par\vskip4pt\def\\{\unskip, }\Coordfont\baselineskip10.5pt\noindent#1}}

 \def\pagecontents{%\TRMargIns new, \Pagetot@l
  \gdef\Pagetot@l{\pagetotal}
  \ifvoid\TRMargIns\else
    \rlap{\kern\hsize\kern10pt\vbox to 0pt{%
         \box\TRMargIns\vss}}\fi
  \ifvoid\topins\else\unvbox\topins\fi
   \dimen@=\dp\@cclv \unvbox\@cclv % open up \box255
   \ifvoid\footins\else % footnote info is present
     \vskip\skip\footins
     \footnoterule
     \unvbox\footins\fi
   \ifr@ggedbottom \kern-\dimen@ \vfil \fi}

  %%%%% Some math accents %%%%%

 \newcount\Ht %pg121; Height register, used in Linefigure & accents

 \def \Acc{\expandafter } %%% What is this for?? WN

 \def\swthat{\raise -1.1 ex\hbox{\sam$\widehat{}$}}
 \def\swttilde{\raise -1.2 ex\hbox{\sam$\widetilde{}$}}
 \def \overdot{{\raise .2 ex \hbox to 0pt {\hss\bf\smash{.}\hss}}}
 \def \overcircle{{\raise .1 ex \hbox to 0pt
    {\sam$\eightpoint\scriptstyle\hss\circ\hss$}}}

 \def \Mathaccent#1#2{{\sam % E.g. #1=\widehat
  \setbox4=\hbox{$\vphantom{#2}$}
  \Ht=\ht4 %pg120
  \setbox5=\hbox{${#1}$}
  \setbox6=\hbox{${#2}$}
  \setbox7=\hbox to .5\wd6{}
  \copy7\kern .1\Ht \raise\Ht sp\hbox{\copy5}\kern-.1\Ht
  \copy7\llap{\box6}
  }}

  \def\SwtCheck #1{
        \ifmmode \check{#1}%
                \else \v {#1}%
                \fi}

 %%  \barpartial : bar over partial is common, tailor!
 \def\barpartial {%
   \kern .17 em
    \overline {\kern -.17 em\partial\kern-.03 em}%
    \kern .03 em}

 %%%   BEtter overline
 
  \def\Overline#1{\setbox1=\hbox{\sam ${#1}$}%
      \ifdim \wd1 > 6pt
    \kern .11 em
    \overline {\kern -.11 em#1\kern-.14 em}
    \kern .14 em
  \else
    \kern .03 em
    \overline {\kern -.03 em#1\kern-.04 em}
    \kern .04 em
  \fi}

 \def\SOverline#1{\setbox1=\hbox{\sam ${#1}$}%
      \ifdim \wd1 > 7pt
    \kern .22 em
    \overline {\kern -.22 em#1\kern-.09 em}%
    \kern .09 em
  \else
    \kern .10 em
    \overline {\kern -.10 em#1\kern-.04 em}%
    \kern .04 em
  \fi}

  %%% Better underline

 \def\Underline#1{\setbox1=\hbox{\sam ${#1}$}%
      \ifdim \wd1 > 6pt
    \kern .11 em
    \underline {\kern -.11 em#1\kern-.14 em}
    \kern .14 em
  \else
    \kern .03 em
    \underline {\kern -.03 em#1\kern-.04 em}
    \kern .04 em
  \fi}

 \def\SUnderline#1{\setbox1=\hbox{\sam ${#1}$}%
      \ifdim \wd1 > 7pt
    \kern .04 em
    \underline {\kern -.04 em#1\kern-.2 em}%
    \kern .2 em
  \else
    \kern .0 em
    \underline {\kern -.0 em#1\kern-.15 em}%
    \kern .15 em
  \fi}

  %%%%% Miscellaneous %%%%%

 \def \Blackbox
   {\leavevmode\hskip .3pt \vbox
   {\hrule height 5pt\hbox{\hskip 4.5pt}}\hskip .5pt}

 \def \XX{\Blackbox\kern.5pt\Blackbox} %% editorial use

  \def\.{.\kern1pt}

  %% unbreakable hyphen (by local change of hyphenchar to -1)
    \def\Hyphen{\edef\this{\the\hyphenchar\font}%
          \hyphenchar\font=-1\char\this\hyphenchar\font=\this}

  %% Prose In Math or Display 
 \ifx\undefined\text
  \def\text#1{\hbox{\rm #1}}\fi %% AMSTeX is more sophisticated

  %% Math Object Names (multi-character math object names)
  %%\nolimits can be cancelled
                                     % by a following \limits if wanted

%%%% Larry's mathsurround hacks:

   \everymath{}  %% initially, but later ...

  \def\PassMath@@{\aftergroup\AfterMath@} %% use \aftergroup LS 5-92

 \let\PassMath@\PassMath@@

 \def\AfterMath@{\futurelet\next\AfterMathMole@}

 \def\AfterMathMole@{%\show\next
      \ifcat\next\space% picks off CR and \par cases too; not \dots
          \def\this{}%{(space)}%
      \else
      \ifcat\next\egroup %
        \def\this{\osumess{Handset mathsurround?? ---(see dollar brace)}}%
      \else
      \def\this{\AAfterMath@}% this minority case slow
      \fi\fi
      \this}

 \def\hyphen@{-}
 \def\paren@{)}
 \def\apostr@{'}

 \def\MSC#1{\kern-.8\mathsurround#1\kern.8\mathsurround}

 \def\AAfterMath@#1{\def\Next{#1}%\show\Next%
    \IN@0\Next @,.;:!?\relax @%
    \ifIN@\def\this{\MSC{\Next}}%
    \else
    \ifx\Next\hyphen@\def\this{\futurelet\next\AfterHyphen@}%
    \else
    \ifx\Next\paren@\def\this{#1}%
    \else 
    \ifx\Next\apostr@\def\this{#1}%
    \else \def\this{\osumess{Handset mathsurround??}%
                 #1}\fi\fi\fi\fi
    \this}

 \def\AfterHyphen@#1{\def\Next{#1}%
   \ifx\Next\hyphen@\def\this{--}\else
   \ifcat\next\space%
   \def\this{\kern-\mathsurround\kern.05em- \Next}\else
   \def\this{\kern-\mathsurround\kern.05em\Hyphen\Next}\fi\fi\this}

%%%% switches
 \def\sam{\mathsurround=\z@\let\PassMath@\relax}  %
 \def\mas{\mathsurround=\StdMathsurround\let\PassMath@\PassMath@@}
 
 \def\Mas{\mathsurround=\StdMathsurround
                \everymath{\PassMath@}\let\PassMath@\PassMath@@}

 \def\m@th{\mathsurround=\z@\everymath{}}%% good general measure

 \def\m@@th{\mathsurround=\z@\everymath={}\let\m@th\relax}

\def\underbar#1{$\setbox\z@\hbox{#1}\dp\z@\z@
      \m@th \underline{\box\z@}$\relax}

\def\mathhexbox#1#2#3{\leavevmode
  \hbox{\m@@th$\m@th \mathchar"#1#2#3$}}

\def\dots{\relax\ifmmode\ldots\else$\m@th\ldots\,$\relax\fi}
   %%% this first \relax is ONLY original

\def\dotfill{\cleaders\hbox{\m@@th$\m@th \mkern1.5mu.\mkern1.5mu$}\hfill}
\def\rightarrowfill{$\m@th\mathord-\mkern-6mu%
  \cleaders\hbox{\m@@th$\mkern-2mu\mathord-\mkern-2mu$}\hfill
  \mkern-6mu\mathord\rightarrow$\relax}
\def\leftarrowfill{$\m@th\mathord\leftarrow\mkern-6mu%
  \cleaders\hbox{\m@@th$\mkern-2mu\mathord-\mkern-2mu$}\hfill
  \mkern-6mu\mathord-$\relax}

\def\downbracefill{$\m@th\braceld\leaders\vrule\hfill\braceru
  \bracelu\leaders\vrule\hfill\bracerd$\relax}
\def\upbracefill{$\m@th\bracelu\leaders\vrule\hfill\bracerd
  \braceld\leaders\vrule\hfill\braceru$\relax}

\def\angle{{\vbox{\m@@th\ialign{$\m@th\scriptstyle##$\crcr
      \not\mathrel{\mkern14mu}\crcr
      \noalign{\nointerlineskip}
      \mkern2.5mu\leaders\hrule height.34pt\hfill\mkern2.5mu\crcr}}}}

\def\big#1{{\m@@th\hbox{$\left#1\vbox to8.5\p@{}\right.\n@space$}}}
\def\Big#1{{\m@@th\hbox{$\left#1\vbox to11.5\p@{}\right.\n@space$}}}
\def\bigg#1{{\m@@th\hbox{$\left#1\vbox to14.5\p@{}\right.\n@space$}}}
\def\Bigg#1{{\m@@th\hbox{$\left#1\vbox to17.5\p@{}\right.\n@space$}}}
\def\n@space{\nulldelimiterspace\z@ \m@th}

\def\root#1\of{\setbox\rootbox\hbox{\m@@th$\m@th\scriptscriptstyle{#1}$}
  \mathpalette\r@@t}
\def\r@@t#1#2{\setbox\z@\hbox{\m@@th$\m@th#1\sqrt{#2}$\relax}
  \dimen@\ht\z@ \advance\dimen@-\dp\z@
  \mkern5mu\raise.6\dimen@\copy\rootbox \mkern-10mu \box\z@}

\def\mathph@nt#1#2{\setbox\z@\hbox{\m@@th$\m@th#1{#2}$}\finph@nt}

\def\mathsm@sh#1#2{\setbox\z@\hbox{\m@@th$\m@th#1{#2}$}\finsm@sh}

\def\@vereq#1#2{\lower.5\p@\vbox{\m@@th\baselineskip\z@skip\lineskip-.5\p@
    \ialign{$\m@th#1\hfil##\hfil$\crcr#2\crcr=\crcr}}}

\def\mathpalette#1#2{\sam\mathchoice{#1\displaystyle{#2}}%
  {#1\textstyle{#2}}{#1\scriptstyle{#2}}{#1\scriptscriptstyle{#2}}\mas}

\def\widehat#1{\setbox\z@\hbox{\sam$#1$}%
 \ifdim\wd\z@>\tw@ em\mathaccent"0\msbfam@5B{#1}%
 \else\mathaccent"0362{#1}\fi}
\def\widetilde#1{\setbox\z@\hbox{\sam$#1$}%
 \ifdim\wd\z@>\tw@ em\mathaccent"0\msbfam@5D{#1}%
 \else\mathaccent"0365{#1}\fi}

 \def\dots{\relax{}
  \ifmmode\def\thedots{\mdots@}\else\def\thedots{\tdots@}\fi %
  \thedots}

 %% \eqno and \leqno need protection
 \let\@ldeqno\eqno\let\@ldleqno\leqno
 \def\eqno{\everymath{}\@ldeqno} \def\leqno{\everymath{}\@ldleqno}

  \let\@ldeqalignno\eqalignno
  \def\eqalignno#1{\sam\@ldeqalignno{#1}\mas}
  \let\@ldeqalign\eqalign
  \def\eqalign#1{\sam\@ldeqalign{#1}\mas}

 \def\overrightarrow#1{\vbox{\m@th\ialign{##\crcr
      \rightarrowfill\crcr\noalign{\kern-\p@\nointerlineskip}
      $\hfil\displaystyle{#1}\hfil$\crcr}}}
 \def\overleftarrow#1{\vbox{\m@th\ialign{##\crcr
      \leftarrowfill\crcr\noalign{\kern-\p@\nointerlineskip}
      $\hfil\displaystyle{#1}\hfil$\crcr}}}
 \def\overbrace#1{\mathop{\vbox{\m@th\ialign{##\crcr\noalign{\kern3\p@}
      \downbracefill\crcr\noalign{\kern3\p@\nointerlineskip}
      $\hfil\displaystyle{#1}\hfil$\crcr}}}\limits}
 \def\underbrace#1{\mathop{\vtop{\m@th\ialign{##\crcr
      $\hfil\displaystyle{#1}\hfil$\crcr\noalign{\kern3\p@\nointerlineskip}
      \upbracefill\crcr\noalign{\kern3\p@}}}}\limits}

  \let\@ldmatrix\matrix
  \let\end@ldmatrix\endmatrix
  \def\matrix{\sam\@ldmatrix}
  \def\endmatrix{\end@ldmatrix\mas}
  \let\@ldgather\gather
  \let\end@ldgather\endgather
  \def\gather{\sam\@ldgather}
  \def\endgather{\end@ldgather\mas}
  \let\@ldalign\align
  \let\end@ldalign\endalign
  \def\align{\sam\@ldalign}
  \def\endalign{\end@ldalign\mas}
  \let\@ldaligned\aligned
  \let\end@ldaligned\endaligned
  \def\aligned{\sam\@ldaligned}
  \def\endaligned{\end@ldaligned\mas}
  \let\@ldtag\tag
  \def\tag{\sam\@ldtag}
   %
  %%% Commutative diagrams : use LamsCD too?

   \let\MinCDArrowWidth\minCDaw@

  %% will be redefined by BoxedEPS.tex

  %%%%% \FigureTitle %%%%%

%%%% End of Larry's mathsurround stuff
%%%% Start of Walter's insert corrections

\newskip\insertskipamount\newskip\inserthardskipamount
\insertskipamount 6pt plus2pt %This is medskipamount without shrink
\inserthardskipamount 6pt
\def\insertskip{\vskip\insertskipamount}
\newcount\SplitTest%        will be set to -1 if a topinsert has split
\def\SetSplitTest{\SplitTest\insertpenalties
  \insert\topins{\floatingpenalty1}%
  \advance\SplitTest-\insertpenalties}
\def\midinsert{\par
 \SaveLastSkip\penalty-150\SetSplitTest\RestoreLastSkip
 \ifnum\SplitTest=-1
  \@midfalse\p@gefalse\else\@midtrue\fi\@ins}
\def\@ins{\par\begingroup\setbox\z@\vbox\bgroup%
  \vglue\inserthardskipamount}
\def\endinsert{\egroup % finish the \vbox
  \if@mid \dimen@\ht\z@ \advance\dimen@\dp\z@
    \advance\dimen@\insertskipamount%            was 12pt (wn)
    \advance\dimen@\pagetotal\advance\dimen@-\pageshrink
    \ifdim\dimen@>\pagegoal\@midfalse\p@gefalse\fi\fi
  \if@mid%
    \ifdim\lastskip<\insertskipamount\removelastskip\insertskip\fi
    \nointerlineskip\box\z@\penalty-200\insertskip
  \else%
    \SaveLastSkip%                                  added (wn)
    \insert\topins{\penalty100 % floating insertion
    \splittopskip\z@skip
    \splitmaxdepth\maxdimen \floatingpenalty\z@
    \ifp@ge \dimen@\dp\z@
    \vbox to\vsize{\unvbox\z@\kern-\dimen@}% depth is zero
    \else \box\z@\nobreak\insertskip\fi}% was \bigskip\fi (wn)
    \RestoreLastSkip%                               added (wn)
   \fi\endgroup}
%% End Walter's insert stuff

 %%%%% Footnotes %%%%%

  \newcount\notenumber
  
  \def\note{\advance\notenumber by 1
    \footnote{\the\notenumber)}}

  \newbox\footbox

 %% The following modifies Plain TeX definitions, qv
  \def\footnote#1{\let\@sf\empty
    %{(the text)} is read later
    \ifhmode\edef\@sf{\spacefactor\the\spacefactor}\/\fi
    \sam${}^{\fam0 #1}$\@sf\vfootnote{#1}}%

  \def\vfootnote#1{\insert\footins\bgroup
     \interlinepenalty100 \splittopskip=1pt
     \floatingpenalty=20000
     \leftskip=0pt\rightskip=0pt%
     \parindent=.3em%% adjust
     \Smallfonts\rm%%osudeG added \Smallfonts
     \FootItem@{#1}%\strut% not nec
     \futurelet\next\fo@t}

  \def\FootItem@#1{\par\hangafter1\hangindent=\FootHang
     \setbox0=\hbox{\ignorespaces#1\unskip}%
     \dimen0=.4em\SetOverhang@% dimen0 is extra space
     \noindent\rlap{\box0}\kern\Overhang\ignorespaces}

  %\MaxFootTag{2)}%% in param file

  \def\fo@t{\ifcat\bgroup\noexpand\next \let\next\f@@t
    \else\let\next\f@t\fi \next}
  \def\f@@t{\bgroup\aftergroup\@foot\let\next}
  \def\f@t#1{\baselineskip=10pt\lineskip=1pt
            \lineskiplimit=0pt #1\@foot}%
     %%osudeG added \baselineskip=? pt\lineskiplimit=0pt
  \def\@foot{%%% special strut osu for end of each note
        \hbox{\vrule height0pt depth5pt width0pt}
        \egroup}
  \skip\footins=12 pt plus 0pt minus 0pt %% was \bigskipamount
    %% space added when footnote is present
  \count\footins=1000 % footnote magnification factor (1 to 1)
  \dimen\footins=8in % maximum footnotes per page

 %%%% Altenatives

  %%  Editorial stuff (delete??)

 \def\osumess#1{\EdSpider{\immediate\write16{Line \the\inputlineno: #1}}}%
 \def\HideEdStuff{\gdef\EdSpider##1{}}

 \font\BigSym=cmmi10 scaled \magstep 4

 \def\change{\InLMargin{\hbox{\BigSym \char63\kern10pt}}}

 \def\beginchange{\InLMargin{\hbox{\sam\twelvepoint$\heartsuit$\kern10pt}}}

 \def\endchange{\InLMargin{\hbox{\sam\twelvepoint$\spadesuit$\kern10pt}}}

 \def\InLMargin#1{\strut\vadjust{%
     \kern-\strutdepth
     \vtop to \strutdepth{%
         \baselineskip\strutdepth
         \llap{\sam$\smash{\hbox{\EdSpider{#1}}}$}\null}}}

 \def\strutdepth{\dp\strutbox}
 \def\strutheight{\ht\strutbox}

 \def\NoteInRMargin#1{\strut\vadjust{%
     \kern-1.001\strutdepth
     \vtop to \strutdepth{%
       \baselineskip\strutdepth
       \vss\rlap{\ninepoint\unskip\hskip\hsize
         \vtop to 0pt{%
           \hsize=16em\hfuzz=\hsize
           \leftskip=10pt%
           \rightskip=0pt plus 10000pt%
           \baselineskip=9.8pt\lineskip=.2pt%
           \let\\\break
           \noindent\EdSpider{#1}\vss}%
                \kern10pt}\hbox{}}%%\hbox{}=\null crucial!!
       }}

 \def\ednote#1{\NoteInRMargin{\tentt #1}}

 \def\cbar{\InLMargin{%
      \dimen0=\strutdepth\advance\dimen0 by \lineskip
      \vrule width 3pt
      height \strutheight depth \dimen0 \kern
      3pt}}

 \def\ccbar{\InLMargin{%
      \dimen0=2\strutdepth\advance\dimen0 by 2\lineskip
      \vrule width 3pt
        height 3\strutheight depth \dimen0 \kern
      3pt}}

 \newinsert\TRMargIns
 \dimen\TRMargIns=\maxdimen
 %\count\TRMargIns=0
 %\skip\TRMargIns=0pt

  \def\Ednote#1{\insert\TRMargIns{%
       \vbox to 0pt{\hsize=140pt\hfuzz=\hsize
           \leftskip=6pt%
           \rightskip=0pt plus 10000pt%
           \baselineskip=9.8pt\lineskip=.2pt%
           \let\\\break
           %\vglue\pagetotal% misplaces notes if inserts are present
           \SetPageRemainder% This ...
           \vglue540pt\vglue-\PageRemainder%  .. is a fix (WN)
           \noindent\EdSpider{\tentt #1}\vss}%
       \smallskip}}

 \def\KillEdStuff{\def\ednote##1{}\def\Ednote##1{}%
      \let\change\relax\let\beginchange\relax\let\endchange\relax
       \let\cbar\relax\let\ccbar\relax}

 %%% Compatibility with osumrip.sty
  %%

 %%% Parameters
  \topskip=12pt
  \newskip\StdBaselineskip % to set \baselineskip
  \StdBaselineskip 12pt
  \lineskip=1.1pt
  \lineskiplimit=.8pt
  \widowpenalty=10000 % 8000 to 10000
  \clubpenalty=10000  % 8000 to 10000
  \abovedisplayskip=6pt plus 1pt minus 1pt
  \abovedisplayshortskip=3pt plus 1.5pt
  \belowdisplayskip=6pt plus 1pt minus 1pt
  \belowdisplayshortskip=5pt plus 1pt minus 1pt
  \hfuzz=1.5pt   % Enable overfull box warnings at console

  \def\StdPretolerance{100}
  \tolerance=\StdPretolerance

  \newdimen\StdMathsurround
  \StdMathsurround=1.5pt % 1pt usual without \Mas
  \mathsurround=\StdMathsurround
  \Mas                   %% sophisticated mathsurround on
 % \Sam                   %% sophisticated mathsurround off

%% marker before English punctuation in displayed math
   \def\prose{\relax\hbox{\kern.6\StdMathsurround}}
  
  \def\StdParskip{0pt}    %% Larry wants {2pt plus 1pt}
  \parskip=\StdParskip
  \parindent=0.5cm
 
%%%% load Times for main body font

  \def\Times{ptmr  } 
  \def\TimesI{ptmri  } 
  \def\TimesB{ptmb  }
  \def\TimesBI{ptmbi  }
  \def\HelveticaN{phvrrn }

  =\Times at 10bp% roman text
  =\TimesB at 10bp% boldface extended
   % slanted roman
  \font\tenit=\TimesI at 10bp% text italic
  =\TimesBI at 10bp

  \font\tenmrm=cmr10  %%new name for math role at full size

%%%%% Fonts at ninepoint %%%%%

    =\Times at 9bp 
    \font\nineit=\TimesI at 9bp 
    =\TimesB at 9bp 
    =\TimesBI at 9bp 

    =\HelveticaN at 9bp 
       % see below

%%%%% Fonts at twelvepoint %%%%%

  =\Times at 12bp
  \font\twelveit=\TimesI at 12bp
  =\TimesB at 12bp

%%%%% Fonts at titlepoint %%%%%

  \font\titleit=\TimesI at 14.4bp
  =\TimesB at 14.4bp

 \SetAuthorHead{AuthorHead} % needs \ninepoint since box set
 \SetTitleHead{TitleHead}  % notably \HeaderFont

%%%% Char adjustments %%%%

  \def\lBr{\raise.125ex\hbox{[\kern.1125ex}}
  \def\rBr{\raise.125ex\hbox{\kern.1125ex]}}

 \setbox\footbox=\hbox{\Smallfonts 2)~}

%% Some optional font dimension and spacing 
%% adjustments beyond this point

%% Correct the lousy spacing of italic f (a hack).

  \bgroup
  \catcode`\@=11 %localised
  \gdef\itSpacing{%
     \xspaceskip=.31em plus.1em minus.05em \sfcode `f=2001
     \itWarning@\let\itWarning@\itWarning@@}
  \gdef\itSpacingOff{%
     \xspaceskip=0pt \sfcode `f=1000
     \let\itWarning@\relax}
   \global\let\itWarning@\relax
  \gdef\itWarning@@{\errmessage{%
  Special italic spacing already in force
  (you have probably omitted an ``endth'').
  See itSpacing macro in osuPSfnt.sty
         }}
  \egroup

 %%% Provisional fontdimen settings
  %%
 \fontdimen1\titlebf=0.0pt
 \fontdimen2\titlebf=3.6135pt
 \fontdimen3\titlebf=2.8908pt
 \fontdimen4\titlebf=1.44539pt
 \fontdimen5\titlebf=6.64882pt
 \fontdimen6\titlebf=14.45398pt
 \fontdimen7\titlebf=1.60439pt

 \fontdimen1\tenbi=0.26794pt
 \fontdimen2\tenbi=2.50937pt
 \fontdimen3\tenbi=2.00749pt
 \fontdimen4\tenbi=1.00374pt
 \fontdimen5\tenbi=4.59717pt
 \fontdimen6\tenbi=10.03749pt
 \fontdimen7\tenbi=1.11415pt

 \fontdimen1\twelverm=0.0pt
 \fontdimen2\twelverm=3.01125pt
 \fontdimen3\twelverm=2.409pt
 \fontdimen4\twelverm=1.2045pt
 \fontdimen5\twelverm=5.39615pt
 \fontdimen6\twelverm=12.045pt
 \fontdimen7\twelverm=1.33699pt

 \fontdimen1\twelveit=0.27731pt
 \fontdimen2\twelveit=3.01125pt
 \fontdimen3\twelveit=2.409pt
 \fontdimen4\twelveit=1.2045pt
 \fontdimen5\twelveit=5.37207pt
 \fontdimen6\twelveit=12.045pt
 \fontdimen7\twelveit=1.33699pt

 \fontdimen1\twelvebf=0.0pt
 \fontdimen2\twelvebf=3.01125pt
 \fontdimen3\twelvebf=2.409pt
 \fontdimen4\twelvebf=1.2045pt
 \fontdimen5\twelvebf=5.5407pt
 \fontdimen6\twelvebf=12.045pt
 \fontdimen7\twelvebf=1.33699pt

 \fontdimen1\tenrm=0.0pt
 \fontdimen2\tenrm=2.50937pt
 \fontdimen3\tenrm=2.00749pt
 \fontdimen4\tenrm=1.00374pt
 \fontdimen5\tenrm=4.49678pt
 \fontdimen6\tenrm=10.03749pt
 \fontdimen7\tenrm=1.11415pt

 \fontdimen1\tenit=0.27731pt
 \fontdimen2\tenit=2.50937pt
 \fontdimen3\tenit=2.00749pt
 \fontdimen4\tenit=1.00374pt
 \fontdimen5\tenit=4.47672pt
 \fontdimen6\tenit=10.03749pt
 \fontdimen7\tenit=1.11415pt

 \fontdimen1\tenbf=0.0pt
 \fontdimen2\tenbf=2.50937pt
 \fontdimen3\tenbf=2.00749pt
 \fontdimen4\tenbf=1.00374pt
 \fontdimen5\tenbf=4.61723pt
 \fontdimen6\tenbf=10.03749pt
 \fontdimen7\tenbf=1.11415pt

 \fontdimen1\ninerm=0.0pt
 \fontdimen2\ninerm=2.25842pt
 \fontdimen3\ninerm=1.80673pt
 \fontdimen4\ninerm=0.90337pt
 \fontdimen5\ninerm=4.0471pt
 \fontdimen6\ninerm=9.03374pt
 \fontdimen7\ninerm=1.00273pt

 \fontdimen1\nineit=0.27731pt
 \fontdimen2\nineit=2.25842pt
 \fontdimen3\nineit=1.80673pt
 \fontdimen4\nineit=0.90337pt
 \fontdimen5\nineit=4.02904pt
 \fontdimen6\nineit=9.03374pt
 \fontdimen7\nineit=1.00273pt

 \fontdimen1\ninebf=0.0pt
 \fontdimen2\ninebf=2.25842pt
 \fontdimen3\ninebf=1.80673pt
 \fontdimen4\ninebf=0.90337pt
 \fontdimen5\ninebf=4.15552pt
 \fontdimen6\ninebf=9.03374pt
 \fontdimen7\ninebf=1.00273pt

 %%% \SetExtraSpaces \MaxSpaceFactor \SetSpaceFactors
  %%  See TeXbook, page 76.

 \newcount\MaxSpaceFactor
 \MaxSpaceFactor=3000 %% to reset later

 %%%%% Tag styles and (hang-) indents
 \def\ItemStyle{\rm}
 \def\NrStyle{\rm}
 \def\ItemItemStyle{\rm}

 %% Analog dimensioning, convenient for local modifications:
 \MaxItemTag{(iii)}
 \MaxItemItemTag{(iii)}
 \MaxNrTag{(2)}
 \MaxFootTag{2)}
 % \MaxReferenceTag{AaaAA} % for biblio
 \def\ReferenceHang{30pt}

 \catcode`\@=\active

%%%%% End of hack of Neumann-Siebenmann macros

\loadbold

=\Times  
=\Times scaled750
=\Times scaled650
\font\rms=\Times scaled 920 

=\TimesBI scaled 860
=\TimesI scaled 860

\textfont0=\rrm  
\scriptfont0=\erm 
\scriptscriptfont0=\srm

\def\Augment#1#2{%
    \toks0\expandafter{#1}\toks2{#2}%
    \edef#1{\the\toks0\the\toks2}}

 \font\twelverma=\Times  scaled 1200
 \font\tenrma=\Times  scaled 1000
 \font\ninerma=\Times scaled 920
 =\Times scaled 840
 \font\sevenrma=\Times scaled 760
 =\Times scaled 680
 \font\fiverma=\Times scaled 600

 \Augment\tenpoint{%
  \textfont0=\tenrma  \scriptfont0=\sevenrma  
  \scriptscriptfont0=\fiverma  }

 \Augment\ninepoint{%
  \textfont0=\ninerma  \scriptfont0=\sevenrma 
  \scriptscriptfont0=\fiverma}

 \Augment\twelvepoint{%
  \textfont0=\twelverma  \scriptfont0=\ninerma  
  \scriptscriptfont0=\sevenrma}

\mathsurround=1pt
\hsize=13.45truecm
\vsize=19.5truecm
\hoffset=1.25truecm
\voffset=2truecm
\advance\baselineskip by 2pt

\predefine\til{\~}
\def\~#1{\relax\ifmmode\widetilde{#1}\else\til{#1}\fi}

\redefine \le{\leqslant}
\redefine \ge{\geqslant}
\define \wt#1{\mathaccent"0365{#1}}
\define \wh#1{\mathaccent"0362{#1}}

\define \iss{\,\Mathaccent{\raise -.8 ex\hbox{$\widetilde{}$\kern.1em}}\rightarrow\,}

\define \ur{\mathop{\fam0 ur}}

\define \minn{\operatorname{\fam0 min\,}}

\define \tpp{\mathop{\fam0 top}}

\define \innf{\operatorname{\fam0 inf\,}}

\define \Tr{\operatorname{\fam0 Tr\,}}

\define \Gal{\mathop{\fam0 Gal}}

\define \Sw{\operatorname{\fam0 Sw}}
\define \sw{\operatorname{\fam0 sw}}
\define \ksw{\operatorname{\fam0 ksw}}

\define \Fil{\operatorname{\fam0 fil}}

\define \ins{\mathop{\fam0 ins}}

\Mas
\HideEdStuff
\rm 
 
%%%% For GT headers and footers:

\def\issn{{\nineit ISSN 1464-8997 (on line) 1464-8989 (printed)}}

\def\gtp{{\nineit Published 10 December 2000: \ \copyright\ Geometry \& 
Topology Publications}}

\def\gtv3{{\nineit Geometry \& Topology Monographs, Volume 3 (2000) --
Invitation to higher local fields}}

%%%%% For section idents:

\def\lione
{{\rms Geometry \& Topology Monographs}}

\def \litwo{{\rms Volume 3: Invitation to higher local fields
}} 

\def\tinfo #1.#2.#3-#4
{{
\noindent  {\lione} \hfill 
\par 
\vskip-1.5pt
\noindent {\litwo} \hfill
\par 
\vskip-1,5pt
\noindent {\rms Part #1, section #2, pages #3--#4} \hfill
\vskip24pt 
}}

\def\tinfos #1.#2.#3-#4
{{
\noindent  {\lione} \hfill 
\par 
\vskip-1.5pt
\noindent {\litwo} \hfill
\par 
\vskip-1.5pt
\noindent {\rms Pages #3--#4} \hfill
\vskip24pt 
}}

\def\tinfoi #1
{{
\noindent  {\lione} \hfill 
\par 
\vskip-1.5pt
\noindent {\litwo} \hfill
\par 
\vskip-1.5pt
\noindent {\rms Pages iii--xi: Introduction and contents} \hfill
\vskip26pt 
}}

%%%% Set headers and footers %%%%

  \def\titlepagehead{\hfil}

  \newif\iftitlepage\titlepagefalse
  \newif\ifblankpage\blankpagefalse
  \def\makeheadline{
     \ifblankpage{}\else%
     \iftitlepage
\vbox{\line{\vbox to 8.5pt{}
\ninerm
\copy\HLinebox \hfill
\hglue5mm\ninebf\folio 
\titlepagehead}}%
      \else
\vbox{\ifodd\pageno\rightheadline\else\leftheadline\fi}%
      \fi\vskip 12pt\fi}%
     \def\rightheadline{\line{\vbox to 8.5pt{}%
      \ninerm
\copy\TitleBox \hfill
\hglue5mm\ninebf\folio}}%
     \def\leftheadline{\line{\vbox to 8.5pt{}%
        \unskip\ninerm\unskip\ninebf\folio\hglue5mm
      %*%
 \hfill \copy\AuthorBox
%\hfill
}}

 \footline={\ifblankpage{}\else
\iftitlepage\ninepoint\sam\hfill%} 
\line{\vbox to 8.5pt{}%\ninerm
\copy\TFLinebox
\hfill
\hglue5mm %\ninebf\folio
}
            \else
\ninepoint\sam\hfill%}
\line{\vbox to 8.5pt{}%\ninerm
\copy\FLinebox
\hfill 
\hglue5mm
}
\hfil\fi\global\titlepagefalse\fi}

\def\blankpage{{\blankpagetrue\noindent\hbox to 10pt{\hss}\vfill
\pagebreak}}

\tenpoint\rm %% always start here
 
  %%% all done and macros loaded!

\pageno=151

\tinfo I.18.151-164

\SetTFLinebox{\gtp }
\SetFLinebox{\gtv3 }
\SetHLinebox{\issn}

\H 18. On ramification theory of monogenic extensions

Luca Spriano

\SetAuthorHead{L. Spriano}
\SetTitleHead{Part I. Section 18. On ramification theory 
 of monogenic extensions \qquad\qquad}

\phantom{}
\par

We discuss ramification theory for finite extensions $L/K$ of a 
complete discrete valuation field $K$.
This theory deals with 
quantities which measure wildness of ramification, such as different, 
the Artin 
(resp.\ Swan) characters and the Artin (resp.\ Swan) conductors.
 When the residue field extension
$k_L/k_K$ is separable there is a complete theory, e.g.\ \cite{S},
but  in general it is not so. In the classical case
(i.e.\ $k_L/k_K$ separable) proofs of many results in ramification theory use the property  that all finite extensions of valuation rings $\Cal O_L/\Cal O_K$ are 
monogenic which is not the case in general. Examples (e.g.\ \cite{Sp})
show that the classical theorems do not hold in general.
Waiting for a beautiful and general ramification theory, we 
consider a class of  extensions $L/K$ which has a 
good ramification theory. We describe this class and we will call its elements 
{\it well ramified extensions}.  All classical results are generalizable for 
well ramified extensions, for example
a generalization of the Hasse--Arf theorem proved by J.~Borger.
We  also concentrate
our attention on other ramification invariants, more 
appropriate and general;
in particular,  we consider two ramification invariants:
the Kato conductor and Hyodo depth of ramification. 

Here we comment on some works on general ramification theory.

The first direction aims to generalize classical ramification invariants
to the general case working with (one dimensional) rational valued 
invariants. 
In his papers de Smit gives some properties about 
ramification jumps and  considers the different and differential
\cite{Sm2}; he
generalizes  the Hilbert formula by using the monogenic conductor \cite{Sm1}.
 We discuss  works of Kato \cite{K3-4} in subsection
18.2.
In \cite{K2} Kato describes   
ramification theory for two-dimensional local fields and
he proves an analogue of the  Hasse--Arf theorem for those Galois extensions
in which the extension of the valuation rings (with respect to the 
discrete valuation of rank 2) is monogenic.

The second direction aims to extend ramification invariants
from one dimensional to either higher dimensional 
 or to more complicated objects
which involve differential forms (as in Kato's works \cite{K4}, \cite{K5}).  
By using higher
local class field theory, Hyodo \cite{H} defines generalized ramification invariants,
like depth of ramification 
(see Theorem 5 below). 
We discuss relations of his invariants with the (one dimensional) Kato conductor in subsection 18.3 below. 
Zhukov \cite{Z} generalizes  the  classical ramification 
theory to the case where $|k_K:k_K^p|=p$ 
(see section 17 of this volume). 
From the viewpoint of this section the existence of Zhukov's theory 
is in particular due to the fact that in  the case
where $|k_K:k_K^p|=p$ one can reduce various assertions to the well ramified case.

\HH 18.0.   Notations and definitions

In this section we recall some general definitions. We only consider
complete discrete valuation fields $K$
with residue fields $k_K$  of characteristic $p>0.$ We also assume that 
$|k_K:k_K^p|$ is finite.

\df Definition

Let $L/K$ be a finite Galois extension, $G=\Gal(L/K)$.
Let 

\noindent $G_0=\Gal(L/L\cap K_{\ur})$ be the inertia subgroup of $G$.
Define functions
$$i_G,s_G\colon G\to \Bbb Z$$ by
$$i_G(\sigma)=\cases
\displaystyle\innf_{x\in \Cal O_L\setminus\{0\}} v_L(\sigma(x)-x)\quad&
\text{if $\sigma\not=1$}\\
+\infty\quad &\text{if $\sigma=1$}
\endcases
$$
and 
$$s_G(\sigma)=\cases
\displaystyle\innf_{x\in \Cal O_L\setminus\{0\}} v_L(\sigma(x)/x-1)\quad&
\text{if $\sigma\not=1$, $\sigma\in G_0$}\\
+\infty\quad &\text{if $\sigma=1$}\\
0\quad &\text{if $\sigma\not\in G_0$}.
\endcases
$$
\enddf
Then $s_G(\sigma)\le i_G(\sigma)\le s_G(\sigma)+1$
and if $k_L/k_K$ is separable, then
$i_G(\sigma)= s_G(\sigma)+1$ for $\sigma\in G_0$. 
Note that the functions $i_G,s_G$ depend not only on the group $G$,
but on the extension $L/K$; 
we will denote $i_G$  also by $i_{L/K}$.

\df Definition

The Swan function  is defined as 
$$
\Sw_G(\sigma)=
\cases
-|k_L:k_K|s_G(\sigma),\quad &\text{if $\sigma\in G_0\setminus\{1\}$}\\
-\displaystyle\sum_{\tau\in G_0\setminus\{1\}}\Sw_G(\tau) ,\quad &
\text{if $\sigma=1$}\\
0\quad &\text{if $\sigma\not\in G_0$}.
\endcases
$$
For a character $\chi$ of $G$
its Swan conductor 
$$
\sw(\chi)=\sw_G(\chi)=(\Sw_G,\chi)=\frac{1}{|G|} \sum_{\sigma\in G}\Sw_G(\sigma)\chi
(\sigma)\leqno(1)$$
is an integer if 
$k_L/k_K$ is separable (Artin's Theorem)
and is not an integer in general (e.g.\  \cite{Sp, Ch. I}).
\enddf

\HH 18.1.  Well ramified extensions

\df Definition

Let $L/K$ be a finite Galois $p$-extension.
The extension $L/K$ is called {\it well ramified}
if $\Cal O_L=\Cal O_K[\alpha]$ for some $\alpha \in L$.
\enddf

\HHH 18.1.1. Structure theorem for well ramified extensions

\df Definition

We say that an extension $L/K$ is in {\it case~I}  if 
  $k_L/k_K$ is separable;
an extension $L/K$ is in  {\it case~II}  
if $|L:K|=|k_L:k_K|$ (i.e.\ $L/K$ is ferociously ramified
in the terminology of 17.0) and $k_L=k_K(a)$ is purely inseparable over $k_K$.
\enddf 

Extensions in case~I and case~II are well ramified.
An extension which is simultaneously in case~I and case~II is the trivial extension. 

We characterize 
well ramified extensions by means of the function $i_G$ in the following
theorem.

\th Theorem 1 {{\rm(\cite{Sp, Prop. 1.5.2})}}

Let $L/K$ be a finite Galois $p$-extension.
Then the following properties are equivalent:
\Roster 
\Item{(i)} $L/K$ is well ramified;

\Item{(ii)} for every normal subgroup $H$ of $G$
the Herbrand property holds:
 
\Item{} for every $1\not=\tau\in G/H$
$$i_{G/H}(\tau)=\dsize \frac{1}{e(L|L^H)} \sum_{\sigma\in \tau H} i_G(\sigma);$$

\Item{(iii)} the Hilbert formula holds:
$$v_L(\Cal D_{L/K})=\sum_{\sigma\not=1}i_G(\sigma)=\sum_{i\ge 0}(|G_i|-1),$$
for the definition of $G_i$ see subsection 18.2. 
\endRoster
\endth

From the definition we immediately deduce that
 if $M/K$ is a Galois subextension of a well ramified $L/K$ then $L/M$ is well ramified; from 
(ii) we conclude that  $M/K$  is well ramified.

Now we consider  well ramified extensions $L/K$ which are not in case~I nor 
in case~II.

\eg Example 

(Well ramified extension not in case~I and not in case~II).
Let $K$ be a complete discrete valuation field of characteristic zero.
 Let $\zeta_{p^2}\in K.$ Consider 
a cyclic extension of degree $p^2$ defined by $L=K(x)$ where $x$ a root 
of the  polynomial
$f(X)=X^{p^2}-(1+u\pi) \alpha^p,\  \alpha\in U_K,\ \Overline{\alpha}\not\in k_K^p,\ \ 
u\in U_K,$ $\pi$ is a prime of $K$. 
Then $e(L|K)=p=f(L|K)^{\ins},$ so $L/K$ is not in case~I nor
in case~II.
Using Theorem 1, one can  show that $\Cal O_L= \Cal O_K[x]$ by checking the 
  Herbrand property.
% in particular for this example we have
%$i_G(\sigma)=e_K/(p-1)$ and $i_G(\sigma^p)=pe_K/(p-1).$
\endeg

\df Definition 

A well ramified extension which is not in case~I and is not  
in case~II is said to be in {\it case~III}. 
\enddf

Note that in case~III
we have $e(L|K)\ge p, f(L|K)^{\ins}\ge p$.

\th Lemma 1

If $L/K$ is a well ramified Galois extension% in  case~III%
, then for every  
ferociously ramified Galois subextension  $E/K$ such that $L/E$ is totally ramified
either $E=K$ or $E=L$.
\endth

\pf Proof

Suppose that there exists $ K\not=E\not= L$,  such that $E/K$
is ferociously ramified and $L/E$ is totally ramified.
 Let $\pi_1$ be a prime of $L$ such that $\Cal O_L=
\Cal O_E[\pi_1].$ Let $\alpha\in E$ be such that $\Cal O_E=
\Cal O_K[\alpha].$  Then we have 
$\Cal O_L=\Cal O_K[\alpha,\pi_1].$
Let $\sigma$ be a $K$-automorphism of $E$ and denote $\~\sigma$ a 
lifting of $\sigma$
to $G=\Gal(L/K).$ It is not difficult to show that   
$i_G(\~\sigma)=\min \{v_L(\~\sigma\pi_1-\pi_1),v_L(\sigma\alpha-\alpha)\}.$
We show that $i_G(\~\sigma)=v_L(\~\sigma\pi_1-\pi_1).$
Suppose we had $i_G(\~\sigma)=v_L(\sigma\alpha-\alpha),$ then 
$$\frac{i_G(\~\sigma)}{e(L|E)}=v_E(\sigma\alpha-\alpha)=i_{E/K}(\sigma).\leqno(*)$$
Furthermore, by Herbrand property we have
$$i_{E/K}(\sigma)=\frac{1}{e(L|E)}\sum_{s\in\sigma\Gal(L/E)} i_G(s)=
\frac{i_G(\~\sigma)}{e(L|E)}+\frac{1}{e(L|E)}
\sum_{ s\not=\tilde\sigma}
i_G(s).$$
So  from $(*)$ we deduce that
$$\frac{1}{e(L|E)}\sum_{ s\not=\tilde\sigma}
 i_G(s)=0,$$
but this is not possible because $i_G(s)\ge 1$ for all $s\in G.$ We have 
shown that
$$
i_G(s)=v_L(s \pi_1-\pi_1) \qquad \text{for all}\ \ s\in G.\leqno(**)
$$
Now note that   $\alpha\not\in 
 \Cal O_K[\pi_1].$ Indeed, from $\alpha=\sum a_i \pi_1^i,\ a_i\in
 \Cal O_K,$ we deduce $\alpha\equiv a_0\ 
(\text{\tenrm mod} \pi_1)$ which is impossible. 
By $(**)$ and  the Hilbert formula (cf.\ Theorem 1)  we have 
$$
v_L(\Cal D_{L/K})=\sum_{s\not=1} i_G(s)=
\sum_{s\not=1} v_L(s\pi_1-\pi_1)= v_L(f'(\pi_1)),\leqno(***)
$$
 where
$f(X)$ denotes the minimal polynomial of $\pi_1$ over $K.$

Now let  the ideal 
$\Cal T_{\pi_1}=\{ x\in\Cal O_L\ : \
x\Cal O_K[\pi_1]\subset \Cal O_L \}$
be the conductor of $\Cal O_K[\pi_1]$ in $\Cal O_L$ (cf.\ \cite{S, Ch. III, \S 6}). We have (cf.\ loc.cit.)
$$\Cal T_{\pi_1}\Cal D_{L/K}=f'(\pi_1)\Cal O_L$$
and then $(***)$ implies $\Cal T_{\pi_1}=\Cal O_L$,    
$\Cal O_L=\Cal O_K[\pi_1]$, which contradicts  
$\alpha\not\in\Cal O_K[\pi_1]$. 
\qed
\endpf

\th Theorem 2 {{\rm (Spriano)}}

Let $L/K$ be a Galois well ramified $p$-extension. 
Put $K_0=L\cap K_{\ur}$.
Then there is a Galois 
subextension $T/K_0$ of $L/K_0$  such that $T/K_0$ is in  case~I and $L/T$ in case~II.
\endth

\pf Proof 

Induction on $|L:K_0|$.

Let $M/K_0$ be a Galois subextension of $L/K_0$ such that
$|L:M|=p$.
Let $T/K_0$ be a Galois subextension of $M/K_0$ such that
$T/K_0$ is totally ramified and $M/T$ is in case~II.
Applying Lemma 1 to $L/T$ we deduce that
$L/M$ is ferociously ramified, hence in case~II.
\qed
\endpf

In particular, if $L/K$ is a Galois $p$-extension
in case~III such that $L\cap K_{\ur}=K$, then
there is a Galois subextension $T/K$ of $L/K$ such that
$T/K$ is in  case~I, $L/T$ in case~I and $K\not=T\not=L$.

\HHH 18.1.2.   Modified ramification function for well ramified extensions

\phantom{}\par

In the general case one can define a filtration of ramification 
groups as follows.
Given two integers    $n, m\ge 0$   the $(n,m)$-ramification 
group   $G_{n,m}$ of $L/K$ is 
$$G_{n,m}=\{\sigma \in G:\ v_L(\sigma(x)-x)\ge n+m,\ \text{for all }\ 
x\in \Cal M_L^m\}.$$

Put     $G_n=G_{n+1,0}$
and $H_n=G_{n,1},$ so that the classical  ramification groups  
are the $G_n$.
It is easy to show that $H_i\ge G_i\ge H_{i+1}$ for $i\ge 0$.
 
In case~I we have $G_i=H_i$ for all $i\ge 0;$ in case~II we have
$G_i=H_{i+1}$ for all $i\ge 0,$ see \cite{Sm1}.
If $L/K$ is in case~III, we leave to the reader the proof of
the following equality
$$
\aligned
&G_i=\{\sigma\in \Gal(L/K): v_L(\sigma(x)-x)\ge i+1 \quad
\text{for all $x\in \Cal O_L$}\}\\
&=\{\sigma\in \Gal(L/K): v_L(\sigma(x)-x)\ge i+2 \quad
\text{for all $x\in \Cal M_L$}\}=H_{i+1}.
\endaligned
$$

We introduce another filtration which
 allows us to simultaneously deal with   case I, II and III. 

\df Definition 

Let $L/K$ be a finite Galois well ramified extension. The modified $t$-th 
{\it ramification group} $G[t]$ for $t\ge 0$ is defined by
$$G[t]=\{\sigma\in \Gal(L/K)\ : \ \  i_G(\sigma)\ge  t\}.$$
We call an integer number $m$ a {\it modified ramification jump} of $L/K$ if 
$G[m]\not=G[m+1]$. 
 
\enddf

From now on we will consider only $p$-extensions.
\df Definition

For a well ramified extension $L/K$ define the modified Hasse--Herbrand function
 ${\goth s}_{L/K}(u)$,
$u\in \Bbb R_{\ge 0}$
as $${\goth s}_{L/K}(u)=\int_0^u \frac{|G[t]|}{e(L|K)}\ dt .$$ 

Put $g_i=|\Cal G_i|$.
If $m\le u\le m+1$  where $m$ is a non-negative integer, then
$${\goth s}_{L/K}(u)=\frac{1}{e(L|K)}(g_1+\cdots +g_m+g_{m+1}(u-m)).$$
 We drop the index
$L/K$ in ${\goth s}_{L/K}$ if there is no risk of confusion. One can show
that the function ${\goth s}$ is continuous, piecewise linear, increasing
and convex.  In case~I, if 
$\varphi_{L/K}$ denotes the classical Hasse--Herbrand function as in  
\cite{S, Ch. IV}, then ${\goth s}_{L/K}(u)=1+\varphi_{L/K}(u-1).$
 We define a {\it modified upper numbering} for ramification groups by
$G({\goth s}_{L/K}(u))=G[u]$.

If $m$ is a modified ramification 
jumps, then
 the number ${\goth s}_{L/K}(m)$ is called a {\it  modified upper ramification jump} 
of $L/K$.
\enddf

For well
ramified extensions  we can show the Herbrand theorem as follows.

\th Lemma 2

For $u\ge  0$ we have
${\goth s}_{L/K}(u)=\displaystyle\frac{1}{e(L|K)}
\displaystyle\sum_{\sigma\in G} \inf (i_G(\sigma),u).$

\endth

The proof goes exactly as in \cite{S, Lemme 3, Ch.IV, \S 3}.

\th Lemma 3

Let $H$ be a normal subgroup of $G$ and $\tau\in G/H$ and let $j(\tau)$
be the upper bound of the integers $i_G(\sigma)$ where $\sigma$ runs over all
automorphisms of $G$ which are congruent to $\tau$ modulo $H$. Then we have
$$i_{L^H/K}(\tau)={\goth s}_{L/L^H}(j(\tau)).$$

\endth   

For the proof see Lemme 4 loc.cit.\ (note that Theorem 1 is fundamental 
in the proof). In order to show Herbrand theorem,
we have to show the multiplicativity in the tower of extensions
of the function ${\goth s}_{L/K}$.

\th Lemma 4

With the above notation,  we have  ${\goth s}_{L/K}={\goth s}_{L^H/K}\circ
 {\goth s}_{L/L^H}.$
\endth

For the proof see Prop. 15 loc.cit.
 
\th Corollary

If $L/K$ is well ramified
and $H$ is a normal subgroup of $G=\Gal(L/K)$, then the 
Herbrand theorem holds:
$$(G/H)(u)=G(u)H/H \quad \text{for all $u\ge 0$}.$$
\endth

It is known that the upper ramification jumps (with respect the classical 
function $\varphi$)  of an abelian   extension
in case~I are integers. 
This is the Hasse--Arf theorem. Clearly the same result holds 
with respect the function ${\goth s}.$ In fact, if  $m$ is  a classical 
ramification jump and  $\varphi_{L/K}(m)$ is the upper ramification jump, then
the modified ramification jump is $m+1$ and the modified upper ramification jumps is 
${\goth s}_{L/K}(m+1)=1+\varphi_{L/K}(m)$ which is an integer. In case II it is 
obvious that the modified upper ramification jumps are integers. 
  For well ramified extensions we have
the following theorem, for the proof see the end of  18.2.

\th Theorem 3 {{\rm (Borger)}}

The modified upper ramification jumps  of  abelian well ramified
extensions are integers.
\endth

\HH 18.2.  The Kato conductor

 We have already remarked that the  Swan conductor $\sw(\chi)$ for
a character $\chi$ of the Galois group $G_K$ is not an integer in general.
In \cite{K3} Kato defined 
a modified Swan conductor  in case~I, II  for any character $\chi$ of $G_K$;  
and  \cite{K4} contains a definition of  an integer valued  conductor (which we will call
 the Kato conductor) for 
characters of degree 1 in the general case i.e.\ not only in cases~I and II.

 We recall its definition. 
The map $K^*\to H^1(K,\Bbb Z/n(1))$ (cf.\ the definition of $H^q(K)$  in subsection 5.1)
induces
a pairing
$$ \{\,\,,\,\,\}\colon H^q(K)\times K_r(K)\to H^{q+r}(K),$$
which we briefly explain only for $K$ of characteristic zero, in 
characteristic $p>0$ see \cite{K4, (1.3)}.
 For  $a\in K^*$ and a fixed $n\ge 0,$ let 
$\{a\}\in H^1(K,\Bbb Z/n(1))$  be
the image  under the connecting homomorphism $K^*\rightarrow$
 $ H^1(K,\Bbb Z/n(1))$ 
induced by the exact sequence of $G_K$-modules
$$1\longrightarrow \Bbb Z/n(1)\longrightarrow K^*_s \buildrel n \over 
\longrightarrow  K^*_s \longrightarrow 1.$$ 
For $a_1,...,a_r\in K^*$  the symbol
$\{a_1,...,a_r\}\in H^r(K,\Bbb Z/n(r))$
is the cup product $\{a_1\}\cup\{a_2\}\cup\cdots\cup\{a_r\}.$
For $\chi\in H^q(K)$ and  $a_1,\dots ,a_r\in K^*$ 
$\{\chi,a_1,...,a_r\}\in H^{q+r}_n(K) $
is the cup product $\{\chi\}\cup\{a_1\}\cup\cdots\cup\{a_r\}.$ Passing 
to the limit we have the element $\{\chi,a_1,...,a_r\}\in H^{q+r}(K).$

\df Definition  

 Following Kato, we define an increasing filtration 
$\{\Fil_m H^q(K)\}_{m\ge 0}$ of $H^q(K)$ by

{}\par

$$\Fil_m H^q(K)=\{\chi\in H^q(K): \{\chi|_M, U_{m+1,M}\}=0
\quad\text{for every $M$} \}
$$
where $M$ runs through all complete discrete valuation fields 
satisfying  $\Cal O_K\subset \Cal O_M$, 
$\Cal M_M= \Cal M_K\Cal O_M;$ here $\chi|_M$ denotes the image of 
$\chi\in H^q(K)$ in $H^q(M).$
\enddf

Then one can show 
$H^q(K)=\cup_{m\ge0} \Fil_mH^q(K)$ \cite{K4, Lemma (2.2)} which allows us to 
give the following definition. 

\df Definition

For $\chi\in H^q(K)$ the {\it Kato conductor} of $\chi$
is the integer $\ksw(\chi)$ defined by 
$$\ksw(\chi)=\min\{m\ge0:\chi\in \Fil_mH^q(K)\}.$$
\enddf
This integer $\ksw(\chi)$ is a generalization of the classical Swan conductor
as stated in the following proposition.

\th Proposition 1

Let   $\chi\in H^1(K)$ and let  $L/K$ be the corresponding \, finite cyclic
extension  and suppose that $L/K$ is in case~I or II.
Then

\item{(a)} $\ksw(\chi)=\sw(\chi)$ (see formula $(1)$).
\item{(b)}
 Let   $t$ be  the maximal modified ramification jump.
Then
$$
\ksw(\chi)=\left\{\aligned {\goth s}_{L/K}(t)-1 \ \ \ &\ \ \text{case I}\\
{\goth s}_{L/K}(t)\ \ \ & \ \ \text{case II.}
\endaligned\right.$$

\endth

\pf  Proof

(a) See \cite{K4, Prop. (6.8)}. (b) This is a computation left to the reader. \qed
\endpf

We compute the Kato conductor in  case~III.

%Borger's theorem can be deduced from

\th Theorem 4 {{\rm (Spriano)}}

If $L/K$ is a cyclic extension in  case~III
and if $\chi$ is the corresponding element of $H^1(K)$, then
$\ksw(\chi)=\sw(\chi)-1$. If $t$ is the maximal modified ramification jump of $L/K$, 
then $\ksw(\chi)={\goth s}_{L/K}(t)-1.$
\endth

\smallskip

Before the proof we explain how to compute the 
Kato conductor $\ksw(\chi)$ where $\chi\in H^1(K)$.  
Consider the pairing 
$H^1(K)\times K^*\to H^2(K),$ ($q=1=r$). It
 coincides with
the symbol $(\cdot, \cdot)$ defined in \cite{S, Ch. XIV}. In particular, if
$\chi\in H^1(K)$ and $a\in K^*,$ then
$\{\chi, a\}=0$ if and only if the element $a$ is a norm of the extension
$L/K$ corresponding to $\chi$. So
 we have to compute the minimal integer $m$ such that $U_{m+1,M}$
is in the  norm of the cyclic extension corresponding to $\chi|_M$ when $M$ 
 runs through all complete discrete valuation fields 
satisfying $\Cal M_M= \Cal M_K\Cal O_M.$  
The  minimal integer $n$ such that $U_{n+1,K}$ is contained in the norm of $L/K$
is not, in general, the
Kato  conductor (for instance if the residue field of $K$ is
algebraically closed)

 Here is a characterization of the Kato conductor which helps to compute it
and does not involve extensions $M/K$, cf.\ \cite{K4, Prop. (6.5)}.

\th Proposition 2

 Let $K$ be a complete discrete valuation field.
\ Suppose that 
 $|k_K:k_K^p|=p^c<\infty,$ and $H^{c+1}_p(k_K)\not=0.$ Then for $\chi\in
 H^q(K)$ and $n\ge 0$  
$$\chi\in\Fil_n H^q(K)\quad \Longleftrightarrow\quad \{\chi, U_{n+1}K_{c+2-q}^M(K)
\}=0\quad  \text{in}\ \ H^{c+2}(K),$$
for the definition of $U_{n+1}K_{c+2-q}^M(K)$ see subsection~4.2. 
\endth 

In the following we will only consider characters $\chi$ such that 
the corresponding cyclic extensions $L/K$ are $p$-extension, because 
$\ksw(\chi)=0$ for tame characters $\chi$,  cf.\ \cite{K4, Prop. (6.1)}. 
We can compute the Kato conductor in the following manner.

\smallskip
\th Corollary

Let $K$ be as in Proposition 2.
Let $\chi\in H^1(K)$ and assume that  the corresponding cyclic extension 
$L/K$ is a $p$-extension. 
Then the minimal integer $n$ such that 
$$U_{n+1,K}\subset N_{L/K} L^*$$
is the Kato conductor of $\chi.$ 
\endth

\pf Proof

 By the hypothesis (i.e.\ $U_{n+1,K}\subset N_{L/K} L^*$) we have
$\ksw(\chi)\ge n.$
Now $U_{n+1,K}\subset N_{L/K} L^*,$ implies that 
 $U_{n+1}K_{c+1}(K)$ is contained in 
\ the norm group $N_{L/K} K_{c+1} (L).$ By \cite{K1, II, Cor. at p. 659} we have 
that $\{\chi,U_{n+1}K_{c+1}(K)\}=0 $ in $H^{c+2}(K)$ and so    
by Proposition 2 $\ksw(\chi)\le n$. \qed
\endpf 
  
\pf Beginning of the proof of Theorem 4

Let $L/K$ be an extension in case~III and let $\chi\in H^1(K)$ 
be the corresponding character. 
We can assume that $H^{c+1}_p(k_K)\not=0,$ otherwise we consider the extension
$k=\cup_{i\ge 0} k_K(T^{p^{-i}})$ of the residue field $k_K$, 
preserving a $p$-base,  
for which $H^{c+1}_p(k)\not=0$ (see \cite{K3, Lemma (3-9)}). 

So by the above Corollary we have to compute the minimal 
integer $n$ such that $U_{n+1,K}\subset N_{L/K} L^*.$

Let $T/K$ be the totally ramified extension defined by Lemma 1 (here
$T/K$ is uniquely determined because the extension $L/K$ is cyclic).
Denote by $U_{v,L}$ for $v\in\Bbb R, v\ge  0$ the
group $U_{n,L}$ where $n$ is the smallest integer $\ge  v.$  

If $t$ is the maximal modified ramification jump  of $L/K$, then
$$U_{{\goth s}_{L/T}(t)+1,T}\subset N_{L/T} L^*\leqno(1)$$
because $L/T$ is in  case~II and its Kato conductor is ${\goth s}_{L/T}(t)$
by Proposition 1 (b).
Now consider the totally ramified extension $T/K$.
By \cite{S, Ch. V, Cor. 3 \S 6}  we have
$$N_{T/K}(U_{s,T})=U_{{\goth s}_{T/K}(s+1)-1,K}\ \ \ \text{if}\ \ \
\Gal(T/K)_{s}=\{1\}.\leqno(2) $$
 
Let $t'=i_{T/K}(\tau)$ be the maximal modified ramification jump  of $T/K$.
Let $r$ be the maximum of $i_{L/K}(\sigma)$ where  $\sigma$
runs over all representatives of the coset $\tau\Gal(L/T)$. 
By Lemma 3 $t'={\goth s}_{L/T}(r)$.
Note that $r<t$ (we explain it in the next paragraph), so 
$$t'={\goth s}_{L/T}(r)<{\goth s}_{L/T}(t).
\leqno(3)$$

To show that $r<t$ it suffices to show that
for a generator $\rho$ of $\Gal(L/K)$
$$i_{L/K}(\rho^{p^m})>i_{L/K}(\rho^{p^{m-1}})$$
for $|T:K|\le p^m \le |L:K|$.
Write $\Cal O_L=\Cal O_K(a)$ then
$$\rho^{p^m}(a)-a=\rho^{p^{m-1}}(b)-b,\quad b=\sum_{i=0}^{p-1}\rho^{p^{m-1}i}(a).$$
Then $b=pa+\pi^if(a)$ where $\pi$ is a prime element of $L$,
$f(X)\in \Cal O_K[X]$ 
and $i=i_{L/K}(\rho^{p^{m-1}})$.
Hence $i_{L/K}(\rho^{p^{m}})=v_L(\rho^{p^m}(a)-a)\ge \minn(i+v_L(p), 2i)$,
so $i_{L/K}(\rho^{p^{m}})>i_{L/K}(\rho^{p^{m-1}})$, as required.

Now we use the fact that the number ${\goth s}_{L/K}(t)$ is an integer (by 
Borger's Theorem). We shall show that $U_{{\goth s}_{L/K}(t),K}\subset
 N_{L/K}L^*$. 
 
By $(3)$ we have $\Gal(T/K)_{{\goth s}_{L/T}(t)}=\{1\}$ and so we can apply
 $(2).$
By $(1)$
we have $U_{{\goth s}_{L/T}(t)+1,T}\subset N_{L/T}L^*,$ and by applying the norm
map $N_{T/K}$ we have (by $(2)$)
$$N_{T/K}(U_{{\goth s}_{L/T}(t)+1,T})=
U_{{\goth s}_{T/K}({\goth s}_{L/T}(t)+2)-1,K}\subset N_{L/K}L^*.$$
Thus it suffices to show that the smallest integer
 $\ge  {\goth s}_{T/K}({\goth s}_{L/T}(t)+2)-1 $ is ${\goth s}_{L/K}(t).$ Indeed
 we have
$${\goth s}_{T/K}({\goth s}_{L/T}(t)+2)-1={\goth s}_{T/K}({\goth s}_{L/T}(t))+\frac{2}{|T:K|}
-1={\goth s}_{L/K}(t)- 1+\frac{2}{p^e}$$
where we have used Lemma 4. 
 By Borger's theorem ${\goth s}_{L/K}(t)$
is an integer and thus  we have
shown that $\ksw(\chi)\le{\goth s}_{L/K}(t)- 1.$

Now we need a lemma which is a  key ingredient to deduce Borger's theorem.

\th Lemma 5

Let $L/K$ be a Galois extension in case~III. 
If 
$k_L=k_K(a^{1/f})$ then 

\noindent $a\in k_K\setminus k_K^p$
 where $f=|L:T|=f(L/K)^{\ins}$. Let $\alpha$ be a lifting 
of $a$ in $K$ and let $M=K(\beta)$
where $\beta^f=\alpha$.

If $\sigma\in \Gal(L/K)$ and $\sigma'\in\Gal(LM/M)$ is such that
$\sigma'|_L=\sigma$ then 
$$i_{LM/M}(\sigma')=e(LM|L)i_{L/K}(\sigma).$$
\endth

\pf Proof

(After J. Borger). Note that the extension $M/K$ is  in case~II 
and $LM/M$ is in case~I, in particular it is totally ramified.
 Let $x\in\Cal O_L$ such that $\Cal O_L=\Cal O_K[x].$ 
One can check that $x^f-\alpha\in \Cal M_L\setminus \Cal M_L^2$. 
Let $g(X)$ be the minimal polynomial of $\beta$ over $K$. 
Then $g(X+x)$ is an Eisenstein polynomial over $L$ (because $g(X+x)\equiv
X^f+x^f-\alpha\equiv X^f \mod \Cal M_L$) and $\beta-x$ is a root 
of $g(X+x)$. So $\beta -x$ is a prime of $LM$ and we have
$$i_{LM/M}(\sigma')=v_{LM}(\sigma'(\beta-x)-(\beta-x))=v_{LM}(\sigma'(x)-x)
=e(LM|L)i_{L/K}(\sigma).$$
\qed
\endpf

\pf Proof of Theorem 3 and Theorem 4

Now we deduce simultaneously the formula for the  Kato conductor in case~III and Borger's theorem. 
We compute the 
classical Artin conductor $A(\chi|_M).$ By  the preceding
lemma we have 
 $$\aligned A(\chi|_M)=& \frac{1}{e(LM|M)}\sum_{\sigma'\in \Gal(LM/M)}
 \chi|_M(\sigma')i_{LM/M}(\sigma')\\
=& \frac{e(LM|L)}{e(LM|M)}\sum_{\sigma'\in \Gal(LM/M)}
 \chi|_M(\sigma') i_{L/K}(\sigma)= \frac{1}{e(L|K)}\sum_{\sigma\in G}
 \chi(\sigma)i_{L/K}(\sigma).
\endaligned\leqno(2)$$
 Since $A(\chi|_M)$ is an integer by
Artin's theorem we deduce that the latter expression is an integer.
Now by the well known arguments one deduces the Hasse--Arf property for
$L/K$.

The above argument also shows that the Swan conductor (=Kato
conductor) 
of $LM/M$ is equal to  $A(\chi|_M)-1,$ which shows that $\ksw(\chi)\ge A(\chi|_M)-1=
\goth s_{L/K}(t)-1,$
so  $\ksw(\chi)=\goth s_{L/K}(t)-1$ and Theorem 4 follows.
\qed 
\endpf 

\HH 18.3. More ramification invariants

\HHH 18.3.1.  Hyodo's depth of ramification

This ramification invariant was  introduced by Hyodo in \cite{H}.
 We are interested
in its link with the Kato conductor.

Let $K$ be an $m$-dimensional local field, $m\ge 1.$
Let $t_1,\dots,t_m$ be a system of local parameters of $K$
and let $\bold v$ be the corresponding valuation.

\df Definition

Let $L/K$ be a finite extension. The {\it depth of ramification} of $L/K$
 is
$$d_K(L/K)=\inf \{ {\bold v}(\Tr_{L/K}(y)/y): y\in L^*\}\in\Bbb Q^m.$$
\enddf
The right hand side expression  exists; and, 
in particular, if $m=1$ then $d_K(L/K)=v_K(\Cal D_{L/K})-(1-v_K(\pi_L)),$ 
see \cite{H}.
The main result about the depth is stated in the following theorem 
(see \cite{H, Th. (1-5)}).

\th Theorem 5 {{\rm (Hyodo)}}

Let $L$ be a finite Galois extension of an $m$-dimensional local field $K$.
For $l\ge 1$ define
$$
{\bold j}(l)={\bold j}_{L/K}(l)=
\cases \max\{{\bold i}: 1\le {\bold i} \in \Bbb Z^m,\quad |\Psi_{L/K}(U_{\bold i}K_m^{\tpp}(K))|
\ge p^l\}
\quad &\text{if it exists}\\
0\quad &\text{otherwise}
\endcases
$$
where $\Psi_{L/K}$ is the reciprocity map{{\rm;}}  
the definition of $U_{\bold i}K_m^{\tpp}(K)$  is given in 17.0.
Then
$$(p-1)\sum_{l\ge 1} {\bold j}(l)/p^l \le d_K(L/K) \le (1-p^{-1})\sum_{l\ge 1} {\bold j}(l).
\leqno(3)$$
Furthermore, these inequalities are the best possible 
{{\tenrm(}}cf.\ \cite{H, Prop. (3-4) and Ex. (3-5)}{{\tenrm)}}. 
\endth

For ${\bold i}\in \Bbb Z^m$, let $G^{\bold i}$ be the image of $U_{\bold i}K_m^{\tpp}(K) $ in 
$\Gal(L/K)$ under the reciprocity map $\Psi_{L/K}$. The numbers ${\bold j}(l)$ are 
called jumping number (by Hyodo) and in the classical case, i.e.\ 
$m=1$, they coincide with 
the upper ramification jumps of $L/K$.  

For local fields (i.e.\ 1-dimensional local fields) one can  show that
the first inequality in (3) is actually an equality. 
Hyodo stated  (\cite{H, p.292})  
{\it ``It seems that we can define nice
 ramification 
groups only when the first equality of (3) holds.''} 

For example, 
if $L/K$ is of degree $p$, then the inequalities in (3) are actually
 equalities and in this case we actually have a nice ramification theory. 
For 
an abelian extension $L/K$  \cite{H, Prop. (3-4)} shows that  
the first equality of (3) holds if at most one diagonal 
component of $E(L/K)$ (for the definition see subsection 1.2) is divisible by~$p$.

Extensions in case~I or II verify the hypothesis of Hyodo's proposition, but
it is not so in case~III. We shall show below that the first equality
does not hold in case~III.

\smallskip
\HHH 18.3.2. The Kato conductor and depth of ramification

\phantom{}\par

Consider an $m$-dimensional local field $K,$ $m\ge 1$.  Proposition 2
of 18.2 
shows (if the first residue field is of characteristic $p>0$)  that for
$\chi\in H^1(K),$ $\chi\in \Fil_n H^1(K)$ if and only if
the induced homomorphism 
$K_m(K)\rightarrow \Bbb Q/\Bbb Z$ annihilates $U_{n+1}K_m(K)$ (cf.\ also
 in \cite{K4, Remark (6.6)}).
This also means that  the Kato conductor of the extension $L/K$
 corresponding to $\chi$ is the $m$-th component  of the 
last ramification jump ${\bold j}(1) $ (recall that 
${\bold j}(1)=\max\{{\bold i}: 1\le {\bold i} \in \Bbb Z^m,\quad |G^{{\bold i}}|
\ge p\}$).

\eg Example 

Let $L/K$ as in Example of 18.1.1 and assume that $K$ is a 2-dimensional
 local
 field with the first residue field of characteristic $p>0$ and let
$\chi\in H^1(K)$ be the corresponding character.   Let ${\bold j}(l)_i$ 
denotes the $i$-th component of ${\bold j}(l)$. Then by Theorem 3 and by the above 
discussion we have 
$$\ksw(\chi)={\bold j}(1)_2=\goth s_{L/K}(pe/(p-1))-1=\frac{(2p-1)e}{p-1}-1.$$
If $T/K$ is the subextension of degree $p,$ we have
$$d_K(T/K)_2=p^{-1}(p-1){\bold j}(2)_2 \ \ \Longrightarrow \ \ {\bold j}(2)_2=\frac{pe}{p-1}-1.
$$
The depth of ramification is easily computed: 
$$d_K(L/K)_2=d_K(T/K)_2+d_K(L/T)_2=\frac{(p-1)}{p}\left(\frac{2pe}{p-1}-1
\right).$$
The left hand side of $(3)$ is $(p-1)({\bold j}(1)/p+{\bold j}(2)/p^2),$ so for the 
second component we have 
$$(p-1)\left(\frac{{\bold j}(1)_2}{p}+\frac{{\bold j}(2)_2}{p^2}\right)=2e-\frac{(p^2-1)}{p^2}
\not=d_K(L/K)_2.$$
Thus, the first 
equality in (3) does not hold for the extension $L/K$.
\endeg

If $K$ is a complete discrete valuation (of rank one) field, then in the well ramified case straightforward  calculations show that  
$$e(L|K)d_K(L/K)=\cases
\sum_{\sigma\not=1}s_G(\sigma)\quad 
&\text{case~I,II}\\
\sum_{\sigma\not=1} s_G(\sigma)-e(L|K)+1 \quad
&\text{case~III}
\endcases
$$
Let  $\chi\in H^1(K)$ and assume that  the corresponding
extension $L/K$ is well ramified.  
Let $t$ denote the last ramification jump 
of $L/K$; then from the previous formula and Theorem 4 we have
$$e(L|K)\ksw(\chi)=\cases
d_L(L/K)+t\quad
&\text{case~I,II}\\
d_L(L/K)+t-1\quad &\text{case~III}\endcases 
$$

%Now consider a cyclic $p$-extension $L/K$. We call ramification jumps 
%the numbers $s_G(\sigma^{p^j}), j\ge 0.$
In general case, we can indicate the following relation between the Kato conductor and
Hyodo's  depth of ramification.

\smallskip
\th Theorem 6 {{\rm (Spriano)}}

Let $\chi\in H^1(K,\Bbb Z/p^n)$  
and let $L/K$ be the corresponding 
cyclic extension. Then
$$\ksw(\chi)\le d_K(L/K)+\frac{t}{e(L|K)}$$ 
where $t$ is the maximal modified ramification jump. 

\endth

\pf Proof 

In  \cite{Sp, Prop. 3.7.3} we show that 
$$\ksw(\chi)\le \left[\frac{1}{e(L|K)}\left(\sum_{\sigma\in G} \Sw_{G}(\sigma)
\chi(\sigma)-M_{L/K}\right)\right],\leqno(*)$$
where $[x]$ indicates the integer part of $x\in \Bbb Q$ and 
the integer $M_{L/K}$ is defined by 
$$d_L(L/K)+M_{L/K}=\sum_{\sigma\not=1} s_G(\sigma).\leqno(**)$$ 
Thus, the inequality
in the statement follows from $(*)$ and $(**).$ \qed
\endpf

\Bib        References

\rf{B} J.-L. Brylinski, {Th\'eorie du corps de classes de Kato
et r\^evetements ab\'eliens de surfaces, } Ann. Inst. Fourier, 
Grenoble {33} n.3 (1983), 23--38.

\rf{E} H.P. Epp,  {Eliminating wild ramification,} Invent. Math. 
{19}(1973), 235--249.

\rf{FV} I. Fesenko and S. Vostokov, {Local Fields and Their 
Extensions, } Transl. Math. Monograph A.M.S. 1993

\rf{H} O. Hyodo,  {Wild ramification in the imperfect residue 
field case, } Advanced Stud. in Pure Math. {12}(1987) Galois 
Representation and Arithmetic  Algebraic Geometry, 287--314.

\rf{K1} K. Kato,  {A generalization of local class field 
theory by using K-groups I}, J. Fac. Sci. Univ. Tokyo, Sect. IA {26}(1979), 303-376, II ibid., {27}(1980), 603--683.

\rf{K2} K. Kato,  {Vanishing cycles, ramification 
of  valuations, and class field theory}, Duke Math. J. 55(1987), 629--659.

\rf{K3} K. Kato,  {Swan conductors with differential values}, 
\ Advanced Studies in \ Pure Math. \ {12}(1987), 315--342.

\rf{K4} K. Kato, {Swan conductors for characters of degree one 
in the imperfect residue field case, } \ Contemp.\ Math. {83}(1989), 
101--131.

\rf{K5} K. Kato,  {Class field theory, $D$-modules,
and ramification on higher dimensional schemes, part I},
Amer. J. Math. 116(1994),  757--784.

\rf{Ku1} M. Kurihara, {On two types of discrete valuation
fields}, Comp. Math. {63}(1987), 237--257.

\rf{Ku2} M. Kurihara, {Abelian extensions of an absolutely 
unramified local field with general residue field, } Invent. Math. 
{93}(1988), 451--480.

\rf{Lo} V.G. Lomadze, {On the ramification theory of
two-dimensional local fields, } Math. USSR Sbornik, {37}(1980), 
349--365.

\rf{M} H. Miki, {On $\Bbb Z_p$-extensions of complete
$p$-adic power  series fields and function fields}, J. Fac. Sci. Univ.\
Tokyo, Sect. IA {21}(1974), 377--393.

\rf{S} J.-P. Serre, {Corps Locaux,} Hermann Paris 1967, 3rd 
Edition.

\rf{Sm1} B. de Smit, {Ramification groups of local fields
with imperfect residue class field,} J. Number Theory {44} No.3 
(1993), 229--236.

\rf{Sm2} B. de Smit, {The different and differential of
local fields with imperfect residue field,} Proceedings of Edinburgh Math. 
Soc. {40}(1997), 353--365.

\rf{Sp} L. Spriano, {Well and fiercely ramified extensions of complete 
discrete valuation fields, with applications to the Kato conductor,} Th\`ese
\`a l'Universit\'e Bordeaux I, 1999.

\rf{Z} I.B. Zhukov, {On ramification theory in the imperfect
residue field case,} Preprint, Univ. Nottingham, (1998).

 \endBib

\Coordinates

Luca Spriano
Department of Mathematics \ University of Bordeaux

351, Cours de la Lib\'eration 
33405  Talence Cedex
France

E-mail: Luca.Spriano\@math.u-bordeaux.fr
\endCoordinates 

\vfill
\pagebreak

\end